\documentclass[pdflatex,sn-mathphys-num]{sn-jnl}

\usepackage{graphicx}
\usepackage{amsmath,amssymb,amsfonts}
\usepackage{amsthm}
\usepackage[mathscr]{euscript}
\usepackage{bm}
\usepackage{xcolor}
\usepackage{lmodern}

\IfFileExists{bbm.sty}{\usepackage{bbm}\newcommand{\indic}{\mathbbm{1}}}{\newcommand{\indic}{\mathbf{1}}}

\theoremstyle{plain}
\newtheorem{theorem}{Theorem}[section]
\newtheorem{lemma}[theorem]{Lemma}
\newtheorem{proposition}[theorem]{Proposition}
\newtheorem{corollary}[theorem]{Corollary}
\theoremstyle{definition}

\theoremstyle{remark}
\newtheorem{remark}[theorem]{Remark}

\numberwithin{equation}{section}

\newcommand{\rme}{{\mathrm e}}
\newcommand{\rmi}{{\mathrm i}}
\newcommand{\RRe}{{\mathrm{Re}}}
\newcommand{\rmd}{{\mathrm d}}
\newcommand{\al}{{\alpha}}
\newcommand{\la}{{\lambda}}
\newcommand{\ep}{{\varepsilon}}
\newcommand{\1}{\indic}
\newcommand{\sF}{{\mathcal F}}
\newcommand{\sT}{{\mathcal T}}
\newcommand{\sI}{{\mathcal I}}

\newcommand{\sG}{{\mathcal G}}

\newcommand{\sX}{{\mathcal X}}

\newcommand{\PP}{{\mathbb P}}
\newcommand{\EE}{{\mathbb E}}
\newcommand{\RR}{{\mathbb R}}

\newcommand{\Var}{\operatorname{Var}}
\newcommand{\Cov}{\operatorname{Cov}}
\newcommand{\supp}{\operatorname{supp}}

\newcommand{\orcidid}[1]{{\unboldmath\normalfont\footnotesize(ORCID~#1)}}

\begin{document}

\title[On the fractional Kac--Ornstein--Uhlenbeck process]{On the fractional Kac--Ornstein--Uhlenbeck process}

\author*[1]{\fnm{Oscar} \sur{L\'opez-Alfonso}\,\orcidid{0000-0002-8525-1137}}\email{ojlopeza@unal.edu.co}

\author[2]{\fnm{Nikita} \sur{Ratanov}\,\orcidid{0000-0002-0242-0549}}\email{rtnnkt@gmail.com}

\affil*[1]{\orgdiv{Departamento de Matem\'aticas, Facultad de Ciencias}, \orgname{Universidad Nacional de Colombia}, \orgaddress{\city{Bogot\'a}, \country{Colombia}}}

\affil[2]{\orgname{Chelyabinsk State University}, \orgaddress{\city{Chelyabinsk}, \country{Russia}}}

\abstract{\unboldmath We introduce a \emph{fractional Kac--Ornstein--Uhlenbeck} process
$X^{(\al)}$, $0<\al\le1$, via a Langevin equation with
mean-reversion rate $\beta,$ driven by a telegraph process fractionally integrated of order
$\al$. The solution is given by $x\rme^{-\beta t}$ plus the convolution of a
two-state telegraph signal with a Kummer memory kernel, $K_\al$. For $0<\al<1$ and
$\beta>0$,  in contrast with the case $\al=1$, this kernel changes sign
exactly once and has zero total mass\textup; 
much of what follows is an interpretation of this single fact. 
The process takes values within an explicitly defined band,
which is precisely the support of its distribution. Its second-order structure is given in closed
form, and its stationary version is \emph{antipersistent}\textup: the
spectral density vanishes at the origin as $|\omega|^{2-2\al}$, the
autocovariance has zero integral, and the stationary variance is given in
closed form. The driving process obeys a law of large numbers. Finally, for
symmetric velocities and under Kac's scaling, the process converges for
$\tfrac12<\al\le1$ to the Gaussian process
$x\rme^{-\beta t}+\int_0^tK_\al(t-s)\rmd W(s)$, and we prove that for
$\al\le\tfrac12$ it has no limit in law.}



\maketitle

\section{Introduction}
\setcounter{equation}{0}

\subsection{The problem setting}
The Ornstein--Uhlenbeck process owes its tractability to a single formula: the
solution of the Langevin equation driven by a Brownian motion is
$x\rme^{-\beta t}$ plus the noise filtered through the exponential kernel
$\rme^{-\beta t}$. This paper investigates what remains of that picture if the
Gaussian noise is replaced by a force with a finite propagation speed that ``remembers'' its
own past. More precisely, we study the process $X^{(\al)}$, which is defined by 
\begin{equation}
\label{eq:fKOU}
X^{(\al)}(t)=x-\beta\int_0^tX^{(\al)}(s)\rmd s+\sF^{(\al)}(t),
\qquad
\sF^{(\al)}(t)=\frac{1}{\Gamma(\al)}\int_0^t(t-s)^{\al-1}c_{\ep(s)}\rmd s,
\end{equation}
where $0<\al\le1$, $\beta\ge0$, and $c_{\ep(t)}$ is a velocity switching between
two values at the jump times of a two-state Markov chain. We call $X^{(\al)}$
the \emph{fractional Kac--Ornstein--Uhlenbeck} (fKOU) process.

The solution keeps the shape of the Gaussian one,
\begin{equation}
\label{eq:repr-intro}
X^{(\al)}(t)=x\rme^{-\beta t}+\int_0^tK_\al(t-s)\,c_{\ep(s)}\rmd s,
\end{equation}
but the exponential gives way to a Kummer kernel $K_\al$, and this kernel is
not positive: for $0<\al<1$ and $\beta>0$ it changes sign exactly once and has
zero total mass. Almost everything in this paper is an interpretation of this single fact.

Two edges of the parameter square $(\al,\beta)\in(0,1]\times[0,\infty)$ have been investigated.
When $\al=1,$ the driving term is the standard integrated telegraph process
$\sT(t)=\int_0^tc_{\ep(s)}\rmd s$, the kernel is $K_1(t)=\rme^{-\beta t}$, and
$X^{(1)}$ is the Kac--Ornstein--Uhlenbeck process of \citet{KOU1,KOU2}. When
$\beta=0,$ there is no mean reversion, the kernel takes the form $t^{\al-1}/\Gamma(\al)$, and
$X^{(\al)}-x$ is the fractionally integrated telegraph process $\sF^{(\al)}$
itself, studied by \citet{FCAA}. On both edges the kernel is positive. 
In the interior, where $0<\al<1$ and $\beta>0$, this property fails, and so the interior is not
a simple interpolation between the two edges; Section \ref{sec:results} 
discusses the implications of the kernel changing sign.

The memory in \eqref{eq:fKOU} is of a specific nature. Replacing an ordinary
integral with a fractional one is a common technique used to endow 
a dynamical system with a memory that decays according 
to a power law rather than exponentially. 
The primary tools here are the
Riemann--Liouville integral and the associated Mittag-Leffler functions 
\citep{saxena,MillerRoss}; 
established approaches for describing processes with memory 
also include fractional stochastic differential equations \citep{anh-fract,PirozziMishura} 
and time-changed processes \citep{straka,time-times,LP21,LP}. 
In the framework of time-modification, the memory is generated
by an independent random ``clock'', typically an inverse subordinator, and the
object of study is a Markov process observed at random times.
For the telegraph process this route was taken by \citet{ANG}, who showed that
the time-fractional run-and-tumble equation is solved by the telegraph process
run with the inverse of a stable subordinator\textup; the resulting motion no longer
has finite velocity, and at every time it can be found at any distance with
positive probability. In
\eqref{eq:fKOU}, there is no such random clock: the memory is \emph{deterministic},
entirely carried by the kernel $K_\al$, while the randomness 
is introduced solely through a two-state velocity. 
The phenomena discussed below emerge from the interplay 
between this deterministic memory and the Markovian velocity.

The velocity, in turn, is finite. The ``telegraph process'' owes its place in
probability to Kac's observation \citep{Kac}: the position of a particle moving at speed
$\pm c$ and changing direction at times determined by a Poisson process of rate
$\la$ is described by the telegraph equation; furthermore, in the limit
as $c,\la\to\infty$ (with $c^2/\la\to1$), this position converges to Brownian
motion. The theory was further developed in the monograph by \citet{RK};
related models include telegraph processes in inhomogeneous media \citep{MPRF}, 
randomly accelerated motions \citep{Kelbert}, subordinated versions
\citep{GSA25}, and models with a double-well potential \citep{APT}. The fKOU process
is a finite-velocity counterpart of a Gaussian process with the same
memory; this raises the natural question: does Kac's scaling 
still hold in the presence of a kernel? It
does when $\al>\tfrac12$, the limit being
$x\rme^{-\beta t}+\int_0^tK_\al(t-s)\rmd W(s)$, which for $\beta=0$ is the
Riemann--Liouville fractional Brownian motion; the role of the threshold
$\al=\tfrac12$ is discussed in Section \ref{sec:kac}.

Starting from Section \ref{sec:kernel}, we
assume $\beta>0$ unless otherwise specified; statements that hold for
$\beta=0$ are formulated accordingly.

\subsection{Main results}\label{sec:results}

The representation \eqref{eq:repr-intro} is Theorem \ref{theo:repr}, and it
reduces the study of $X^{(\al)}$ to that of its kernel. Theorem
\ref{theo:kernel} describes the shape of $K_\al$. For $0<\al<1$ and $\beta>0$
the function $\rme^{\beta t}K_\al(t)$ decreases strictly from $+\infty$ to
$-\infty$, so that $K_\al$ has a unique zero $t^*=t^*(\al,\beta)$, is positive
before it and negative after it, and $\int_0^\infty K_\al(s)\rmd s=0$: the
recent past pushes the process in the direction of the velocity, the remote
past pulls it back, and the two effects balance exactly. For $\al=1$ the kernel
is $\rme^{-\beta t}$ and nothing of the kind occurs. The results of this paper
are three readings of this dichotomy.

The first reading concerns the range of the process. For $0<\al<1$,
$X^{(\al)}(t)-x\rme^{-\beta t}$ takes values in an explicit interval
$[m_\al(t),M_\al(t)]$ (Theorem \ref{theo:band}), and this interval is exactly
the topological support of its law, whatever the initial state (Proposition
\ref{prop:support}). As $t\to\infty$ it converges to the symmetric band
$\bigl[-c\,\|K_\al\|_{L^1},\,c\,\|K_\al\|_{L^1}\bigr]$, where $2c=c_0-c_1$: the
limiting band depends on the two velocities only through their difference,
because a kernel of zero mass does not see their mean. For $\al=1$ the band is
$[c_1K_2(t),c_0K_2(t)]$ with $K_2(t)=(1-\rme^{-\beta t})/\beta$, and it
converges to the attractor $[c_1/\beta,c_0/\beta]$ of \citet{KOU1} (Corollary
\ref{cor:attractor}).

The second reading is the second-order theory. The means of $X^{(\al)}(t)$ are
combinations of three Kummer functions and the variances are single integrals
of such functions (Theorems \ref{theo:moments} and \ref{theo:acov}), all within
one two-parameter family of kernels that absorbs both the mean-reversion rate
and the switching rate of the velocity. One consequence has no counterpart at
$\al=1$: the two conditional means cross exactly once, so that starting in the
state of larger velocity eventually lowers the mean (Corollary
\ref{cor:means-cross}). For $0<\al<1$ a stationary version
$\widetilde X^{(\al)}$ exists (Proposition \ref{prop:stationary}), and it is
\emph{antipersistent} (Theorem \ref{theo:antipersistence}): its spectral
density vanishes at the origin like $|\omega|^{2-2\al}$, and its autocovariance
$r$ satisfies $\int_0^\infty r(h)\rmd h=0$. For $\tfrac12<\al<1$ this
low-frequency behaviour is that of a fractional Gaussian noise of Hurst index
$H=\al-\tfrac12$; for $\al\le\tfrac12$ there is no fractional Gaussian noise
to compare with. The stationary variance $r(0)$ is given in closed form
(Proposition \ref{prop:r0}), reducing at $\al=1$ to that of the KOU process.
Finally, the driving process obeys a law of large numbers,
$t^{-\al}\sF^{(\al)}(t)\to\kappa/\Gamma(1+\al)$ almost surely, and so does the
time integral of $X^{(\al)}$ (Proposition \ref{prop:lln}).

The third reading is the scaling limit. Here, and only here, the velocity is
symmetric, $c_0=-c_1=c$ and $\la_0=\la_1=\la$, and we let $c,\la\to\infty$ with
$c^2/\la\to1$. For $\tfrac12<\al\le1$ the process $X^{(\al)}$ converges weakly
in $C[0,T]$ to the Gaussian process
$Y^{(\al)}(t)=x\rme^{-\beta t}+\int_0^tK_\al(t-s)\rmd W(s)$, with $W$ a standard
Brownian motion (Theorem \ref{theo:kac-limit}). The proof is self-contained:
the convergence of the integrated telegraph process is obtained in H\"older
norm from the Cattaneo equation satisfied by its characteristic functions, and
the passage to $X^{(\al)}$ is a continuous mapping, bounded on $\gamma$-H\"older
paths when $\gamma>1-\al$ but not at $\gamma=1-\al$. The threshold
$\al=\tfrac12$ is where the window $1-\al<\gamma<\tfrac12$ closes, and it is
also where $K_\al$ ceases to be square-integrable; below it the variance of
$\sF^{(\al)}$ diverges under the same scaling, and $X^{(\al)}$ has no limit in
law (Proposition \ref{prop:nolimit}).

\subsection{Organization}

The paper is organized as follows. Section \ref{sec:prelim} collects what is
needed about the two-state Markov process, about fractional integration and the
process $\sF^{(\al)}$, and about the special functions in which the formulae
below are written. Section \ref{sec:fKOU} proves that \eqref{eq:fKOU} is well
posed, establishes the convolution representation, and derives from it a
pathwise formula and the regularity of the paths. Section \ref{sec:kernel}
describes the shape of the kernel and deduces the invariant band and its
identification with the topological support. Section \ref{sec:moments} develops
the second-order theory: means and variances in closed form, the autocovariance
at finite time, the stationary version (Section \ref{sec:stationary}) and its
antipersistence (Section \ref{sec:antipersistence}); it closes with a law of
large numbers for $\sF^{(\al)}$ (Section \ref{sec:lln}). Section \ref{sec:kac},
which is independent of Sections \ref{sec:stationary}--\ref{sec:lln}, proves the
scaling limit for $\al>\tfrac12$ and examines the regime $\al\le\tfrac12$.
Section \ref{sec:conclusions} contains some concluding remarks.

\section{Preliminaries}\label{sec:prelim}
\setcounter{equation}{0}

\subsection{The two-state Markov process}

Let $\ep=\ep(t)\in\{0,1\}$, $t\ge0$, be a two-state c\`adl\`ag Markov process
with infinitesimal generator
\begin{equation}
\label{def:Lambda}
\Lambda=\begin{pmatrix}-\la_0&\la_0\\ \la_1&-\la_1\end{pmatrix},
\qquad \la_0,\la_1>0,
\end{equation}
stationary probabilities
$p_0^*=\la_1/(\la_0+\la_1)$, $p_1^*=\la_0/(\la_0+\la_1)$, and transition matrix
$\Pi(t)=(\pi_{ij}(t))$,
\begin{equation}
\label{eq:pi}
\Pi(t)=\rme^{t\Lambda}
=\begin{pmatrix}
p_0^*+p_1^*\rme^{-2\la t}&p_1^*-p_1^*\rme^{-2\la t}\\
p_0^*-p_0^*\rme^{-2\la t}&p_1^*+p_0^*\rme^{-2\la t}
\end{pmatrix},\qquad 2\la:=\la_0+\la_1 .
\end{equation}
We write $\PP_i\{\cdot\}=\PP\{\cdot\mid\ep(0)=i\}$ and $\EE_i[\cdot]$ for the
corresponding expectation.
Let $0=\tau_0<\tau_1<\tau_2<\dots$ be the switching times,
$N(t)=\max\{n\ge0:\tau_n\le t\}$ be the associated counting process, and
\begin{equation}
\label{def:simplex}
\Delta_n(t):=\{(\tau_1,\dots,\tau_n)\in\RR^n:0<\tau_1<\dots<\tau_n<t\}.
\end{equation}

Given $\ep(0)=i$, we write $i_k:=(i+k)\bmod2$ for the state after $k$
switchings, so that $\ep(\tau_k)=i_k$. Four elementary properties of $\ep$ will
be used repeatedly.

\begin{lemma}[Switching structure]\label{lem:switch}
Let $i\in\{0,1\}$.
\begin{enumerate}
\item[\textup{(i)}] \textup{(Switching times).} For $t>0$\textup,
$\PP_i\{N(t)=0\}=\rme^{-\la_it}$\textup, and for $n\ge1$ and every Borel set
$B\subseteq\Delta_n(t)$\textup,
\[
\PP_i\bigl\{N(t)=n,\ (\tau_1,\dots,\tau_n)\in B\bigr\}
=\int_BL_n^{(i)}(s_1,\dots,s_n;t)\,\rmd s_1\cdots\rmd s_n ,
\]
where
\begin{equation}
\label{eq:Lambdan}
L_n^{(i)}(s_1,\dots,s_n;t)=\Bigl(\prod_{k=0}^{n-1}\la_{i_k}\Bigr)
\exp\Bigl(-\sum_{k=1}^{n+1}\la_{i_{k-1}}(s_k-s_{k-1})\Bigr),
\qquad(s_1,\dots,s_n)\in\Delta_n(t);
\end{equation}
where $s_0:=0$ and $s_{n+1}:=t$. In particular\textup, $L_n^{(i)}\le\la_{\max}^n$ and
$\PP_i\{N(t)=n\}\le(\la_{\max}t)^n/n!$\textup, where
$\la_{\max}=\max(\la_0,\la_1)$.
\item[\textup{(ii)}] \textup{(Strong Markov property).}
If $T$ is an almost surely finite stopping time for $\ep$\textup, then
conditionally on $\{\ep(T)=j\}$ the shifted process $\ep(T+\cdot)$ is a copy of
$\ep$ started at $j$\textup, independent of $\sG_T$\textup, where
$\sG_t:=\sigma\bigl(\ep(s):s\le t\bigr)$. The same holds for a pair of
independent copies of $\ep$\textup, which is a Markov jump process on
$\{0,1\}^2$.
\item[\textup{(iii)}] \textup{(Ergodicity).} $\ep$ is positive recurrent with
invariant law $(p_0^*,p_1^*)$\textup, so that for any real $A_0,A_1$\textup,
\begin{equation}
\label{eq:ergodic}
\frac1t\int_0^{t}A_{\ep(s)}\rmd s
\ \longrightarrow\ p_0^*A_0+p_1^*A_1\qquad\text{a.s. as }t\to\infty .
\end{equation}
\item[\textup{(iv)}] \textup{(Stationary version).} $\ep$ admits a two-sided
stationary extension $\{\ep(t)\}_{t\in\RR}$\textup, unique in law\textup, with
$\PP\{\ep(t)=i\}=p_i^*$ for all $t$. It is used in Sections
\ref{sec:stationary} and \ref{sec:antipersistence}\textup, and expectations
under it are denoted by $\EE$.
\end{enumerate}
\end{lemma}
\begin{proof} 
(i) Given $\{\ep(0)=i\},$ the holding times $\tau_{k+1}-\tau_k$ are independent and
exponentially distributed with parameter $\la_{i_k}$, which leads to \eqref{eq:Lambdan}; the final
factor, $\rme^{-\la_{i_n}(t-\tau_n)}$, is the probability that no further switching occurs
up to time $t$. 
The bound follows from
$L_n^{(i)}\le\la_{\max}^n$ and
$|\Delta_n(t)|=t^n/n!$.
Statements (ii) and (iii) are standard for finite-state Markov jump processes;
see e.g.\ \citet[Theorems 13.3.3, 13.4.10 and 13.4.15]{Bremaud}. 
Statement (iv) follows from Kolmogorov's extension theorem applied to the stationary
transition semigroup \eqref{eq:pi}.
\end{proof}

Fix $c_0>c_1$ and set $2c:=c_0-c_1>0,\;
2\bar c:=c_0+c_1,$ and
\begin{equation}
\label{def:constants}
\kappa:=p_0^*c_0+p_1^*c_1,\qquad
\bar\kappa:=p_1^*c_0+p_0^*c_1,\qquad
\sigma^2:=4c^2p_0^*p_1^* .
\end{equation}
The process $t\mapsto c_{\ep(t)}$ is the \emph{simple telegraph process}. Let 
$\eta(t):=1-2\ep(t)\in\{-1, +1\}$; then
\begin{equation}\label{eq:centered-c}
c_{\ep(t)}=\bar c+c\,\eta(t),
\end{equation}
so that $|c_{\ep(t)}|\le M:=\max(|c_0|,|c_1|)=|\bar c|+c$; by
\eqref{eq:ergodic}, $\kappa$ is its stationary mean. Note for later use that
\begin{equation}
\label{eq:kappabar}
\bar\kappa-\kappa=2c\,(p_1^*-p_0^*).
\end{equation}

\begin{lemma}
\label{lem:cov}
For all $t,t_1,t_2>0$,
\begin{equation}
\label{eq:mean01}
\EE_0[c_{\ep(t)}]=\kappa+2cp_1^*\rme^{-2\la t},\qquad
\EE_1[c_{\ep(t)}]=\kappa-2cp_0^*\rme^{-2\la t},
\end{equation}
\begin{align}
\label{eq:cov0}
\EE_0\bigl[c_{\ep(t_1)}c_{\ep(t_2)}\bigr]
&=\kappa^2+2cp_1^*\kappa\rme^{-2\la(t_1\wedge t_2)}
+2cp_1^*\bar\kappa\rme^{-2\la(t_1\vee t_2)}
+4c^2p_0^*p_1^*\rme^{-2\la|t_1-t_2|},\\
\label{eq:cov1}
\EE_1\bigl[c_{\ep(t_1)}c_{\ep(t_2)}\bigr]
&=\kappa^2-2cp_0^*\kappa\rme^{-2\la(t_1\wedge t_2)}
-2cp_0^*\bar\kappa\rme^{-2\la(t_1\vee t_2)}
+4c^2p_0^*p_1^*\rme^{-2\la|t_1-t_2|},
\end{align}
and consequently
\begin{align}
\label{eq:ACov0}
\Cov_0\bigl(c_{\ep(t_1)},c_{\ep(t_2)}\bigr)
&=4c^2p_1^*\Bigl[(p_1^*-p_0^*)\rme^{-2\la(t_1\vee t_2)}
+p_0^*\rme^{-2\la|t_1-t_2|}-p_1^*\rme^{-2\la(t_1+t_2)}\Bigr],\\
\label{eq:ACov1}
\Cov_1\bigl(c_{\ep(t_1)},c_{\ep(t_2)}\bigr)
&=4c^2p_0^*\Bigl[(p_0^*-p_1^*)\rme^{-2\la(t_1\vee t_2)}
+p_1^*\rme^{-2\la|t_1-t_2|}-p_0^*\rme^{-2\la(t_1+t_2)}\Bigr].
\end{align}
Under the two-sided stationary version of Lemma \textup{\ref{lem:switch}},
$\EE[c_{\ep(t)}]=\kappa$ for every $t$ and
\begin{equation}
\label{eq:cov-stat}
\Cov\bigl(c_{\ep(t_1)},c_{\ep(t_2)}\bigr)=\sigma^2\rme^{-2\la|t_1-t_2|},
\qquad t_1,t_2\in\RR .
\end{equation}
\end{lemma}

\begin{proof} To simplify notation, we first take $\ep(0)=0$.

Since
$\EE_0[c_{\ep(t)}]=c_0\pi_{00}(t)+c_1\pi_{01}(t)$, formula \eqref{eq:pi} gives
\[
\EE_0[c_{\ep(t)}]
=c_0\bigl(p_0^*+p_1^*\rme^{-2\la t}\bigr)+c_1\bigl(p_1^*-p_1^*\rme^{-2\la t}\bigr)
=\kappa+p_1^*(c_0-c_1)\rme^{-2\la t},
\]
which is the first identity in \eqref{eq:mean01}; the second is obtained in the
same way from the second row of \eqref{eq:pi}.

All four expressions \eqref{eq:cov0}--\eqref{eq:ACov1} are symmetric in
$(t_1,t_2)$, so we may assume that $0<t_1\le t_2$. 

Conditioning on $\ep(t_1)$ and using the Markov property
together with \eqref{eq:mean01},
\begin{align*}
\EE_0\bigl[c_{\ep(t_1)}c_{\ep(t_2)}\bigr]
&=\sum_{j\in\{0,1\}}c_j\,\pi_{0j}(t_1)\,\EE_j\bigl[c_{\ep(t_2-t_1)}\bigr]\\
&=\kappa\,\EE_0\bigl[c_{\ep(t_1)}\bigr]
+2c\,\rme^{-2\la(t_2-t_1)}
\Bigl[c_0\pi_{00}(t_1)p_1^*-c_1\pi_{01}(t_1)p_0^*\Bigr].
\end{align*}
By \eqref{eq:pi} the expression in brackets equals
\[
c_0p_1^*\bigl(p_0^*+p_1^*\rme^{-2\la t_1}\bigr)
-c_1p_0^*p_1^*\bigl(1-\rme^{-2\la t_1}\bigr)
=2c\,p_0^*p_1^*+p_1^*\bar\kappa\,\rme^{-2\la t_1},
\]
whence, again using \eqref{eq:mean01}, we obtain
\[
\EE_0\bigl[c_{\ep(t_1)}c_{\ep(t_2)}\bigr]
=\kappa^2+2cp_1^*\kappa\rme^{-2\la t_1}+2cp_1^*\bar\kappa\rme^{-2\la t_2}
+4c^2p_0^*p_1^*\rme^{-2\la(t_2-t_1)},
\]
which is \eqref{eq:cov0}.
 Subtracting
$\EE_0[c_{\ep(t_1)}]\EE_0[c_{\ep(t_2)}]
=\kappa^2+2cp_1^*\kappa(\rme^{-2\la t_1}+\rme^{-2\la t_2})
+4c^2(p_1^*)^2\rme^{-2\la(t_1+t_2)}$ from \eqref{eq:cov0} leaves
\[
\Cov_0\bigl(c_{\ep(t_1)},c_{\ep(t_2)}\bigr)
=2cp_1^*(\bar\kappa-\kappa)\rme^{-2\la t_2}
+4c^2p_1^*\Bigl[p_0^*\rme^{-2\la(t_2-t_1)}-p_1^*\rme^{-2\la(t_1+t_2)}\Bigr],
\]
and \eqref{eq:kappabar} turns the first term into
$4c^2p_1^*(p_1^*-p_0^*)\rme^{-2\la t_2}$, giving \eqref{eq:ACov0}.

Interchanging the labels $0$ and $1$ maps
$c_0\mapsto c_1$, $\la_0\mapsto\la_1$, hence $p_0^*\mapsto p_1^*$; it leaves
$\la$, $\kappa$ and $\bar\kappa$ invariant and replaces $c$ by $-c$.
Applying it to \eqref{eq:mean01}, \eqref{eq:cov0} and \eqref{eq:ACov0} yields
the corresponding formulae for $\PP_1$, that is \eqref{eq:cov1} and
\eqref{eq:ACov1}.

Averaging \eqref{eq:mean01} over the initial law
$(p_0^*,p_1^*)$ gives $\EE[c_{\ep(t)}]=\kappa$, and averaging
\eqref{eq:cov0}--\eqref{eq:cov1} in the same way cancels the two terms carrying
$\rme^{-2\la t_1}$ and $\rme^{-2\la t_2}$, leaving
$\EE[c_{\ep(t_1)}c_{\ep(t_2)}]=\kappa^2+\sigma^2\rme^{-2\la|t_1-t_2|}$.
This is \eqref{eq:cov-stat}, and it extends to $t_1,t_2\in\RR$ by stationarity.
\end{proof}

Formulae \eqref{eq:mean01}--\eqref{eq:cov1} are taken from
\citet[Lemma 2.1]{FCAA}; the covariance form
\eqref{eq:ACov0}--\eqref{eq:ACov1} and the stationary identity
\eqref{eq:cov-stat} are the versions we shall actually use, the latter being the
input for the stationary autocovariance in Proposition \ref{prop:acov-stat}.

\subsection{Fractional integration and the process $\sF^{(\al)}$}

For $\al>0$ the left-sided Riemann--Liouville fractional integral of a locally
bounded $g$ is
\begin{equation}
\label{def:Inu}
I^\al[g](t)=\frac{1}{\Gamma(\al)}\int_0^t(t-s)^{\al-1}g(s)\rmd s,
\qquad t>0,
\end{equation}
with $I^0=\mathrm{id}$; for piecewise continuous $g$, which covers $c_{\ep(\cdot)}$,
it satisfies the semigroup property
$I^{\al_1}I^{\al_2}=I^{\al_1+\al_2}$, see \citet[Ch.~III, Theorem~1, and Sect.~IV.6]{MillerRoss}.

The fractionally integrated telegraph process of order $\al\in(0,1]$ is
$\sF^{(\al)}=I^\al[c_{\ep(\cdot)}]$, that is \eqref{eq:fKOU}; it was introduced
and studied by \citet{FCAA}.
For $\al=1$ it is the classical integrated telegraph process
\begin{equation}
\label{def:T}
\sT(t):=\sF^{(1)}(t)=\int_0^tc_{\ep(s)}\rmd s,
\end{equation}
for which we refer to \citet{RK}.

\begin{lemma}[Paths of $\sF^{(\al)}$]\label{lem:pathF}
Let $0<\al\le1$. Almost surely\textup:
\begin{enumerate}
\item[\textup{(i)}] $\sF^{(\al)}$ is continuous\textup, continuously
differentiable on each open interval between consecutive switching
times\textup, and for every $t>0$\textup,
\begin{equation}
\label{eq:pathF}
\sF^{(\al)}(t)=\frac{1}{\Gamma(1+\al)}
\Bigl[c_{\ep(0)}t^\al+2c\sum_{n=1}^{N(t)}(-1)^{\ep(\tau_n)}(t-\tau_n)^\al\Bigr];
\end{equation}
\item[\textup{(ii)}] for $0\le t<t'$\textup,
\begin{equation}\label{eq:holder}
\bigl|\sF^{(\al)}(t')-\sF^{(\al)}(t)\bigr|
\le\frac{2M}{\Gamma(1+\al)}\,(t'-t)^\al;
\end{equation}
\item[\textup{(iii)}] as $t\downarrow0$\textup,
\begin{equation}
\label{eq:Fal0}
t^{-\al}\sF^{(\al)}(t)\longrightarrow\frac{c_{\ep(0)}}{\Gamma(1+\al)}.
\end{equation}
\end{enumerate}
\end{lemma}

\begin{proof}
\textup{(i)} Split the integral defining $\sF^{(\al)}(t)$ at the switching
times. On $[\tau_k,\tau_{k+1})$ the velocity equals $c_{i_k}$, and
\[
\frac1{\Gamma(\al)}\int_{\tau_k}^{\tau_{k+1}\wedge t}(t-u)^{\al-1}\rmd u
=\frac{(t-\tau_k)^\al-(t-\tau_{k+1}\wedge t)^\al}{\Gamma(1+\al)}.
\]
Summation by parts, together with
$c_{i_n}-c_{i_{n-1}}=2c(-1)^{\ep(\tau_n)}$, gives \eqref{eq:pathF}; see also
\citet[\S3]{FCAA}. Each term $(t-\tau_n)^\al$ is continuous on $[\tau_n,\infty)$
and continuously differentiable on $(\tau_n,\infty)$, which is the regularity
asserted.

\textup{(ii)} Splitting the difference at $t$,
\[
\Gamma(\al)\bigl[\sF^{(\al)}(t')-\sF^{(\al)}(t)\bigr]
=\int_0^t\bigl[(t'-s)^{\al-1}-(t-s)^{\al-1}\bigr]c_{\ep(s)}\rmd s
+\int_t^{t'}(t'-s)^{\al-1}c_{\ep(s)}\rmd s .
\]
The second integral is bounded in absolute value by $M(t'-t)^\al/\al$, and
\[
\int_0^t\bigl|(t-s)^{\al-1}-(t'-s)^{\al-1}\bigr|\rmd s
=\frac{t^\al-{t'}^\al+(t'-t)^\al}{\al}\le\frac{(t'-t)^\al}{\al},
\]
which gives \eqref{eq:holder}.

\textup{(iii)} By definition,
\[
\frac{1}{t^\al}\sF^{(\al)}(t)=\frac{1}{\Gamma(\al)}\int_0^t(1-s/t)^{\al-1}c_{\ep(s)}\frac{\rmd s}{t}
=\frac{1}{\Gamma(\al)}\int_0^1(1-u)^{\al-1}c_{\ep(tu)}\rmd u
\longrightarrow\frac{c_{\ep(0)}}{\Gamma(1+\al)},
\qquad t\downarrow0,
\]
by dominated convergence, since $\ep$ is right-continuous at $0$.
\end{proof}

\begin{corollary}[Infinite one-sided derivative]\label{cor:Fcusp}
Let $0<\al<1$. Almost surely\textup, the right derivative of $\sF^{(\al)}$ is
infinite at every switching time\textup, and also at $t=0$ whenever
$c_{\ep(0)}\ne0$.
\end{corollary}

\begin{proof}
By \eqref{eq:pathF}, for $n\ge1$,
\[
\sF^{(\al)}(\tau_n+\delta)-\sF^{(\al)}(\tau_n)=2c(-1)^{\ep(\tau_n)}\frac{\delta^\al}{\Gamma(1+\al)}+O(\delta),
\qquad \delta\downarrow0,
\]
so the right derivative at $\tau_n$ is infinite. The same holds at $t=0$ by
\eqref{eq:Fal0}.
\end{proof}

\begin{remark}[Loss of the Markov property]\label{rem:nonmarkov}
For $0<\al<1$ the value of $\sF^{(\al)}(t)$ depends, by \eqref{eq:pathF}, on the
whole switching history $\{\tau_n\}_{n\le N(t)}$ through the weights
$(t-\tau_n)^{\al}$, and the same is true of $X^{(\al)}$. Neither process is
Markov, and neither is the pair it forms with $\ep$\textup; this is proved in
a companion paper (L\'opez-Alfonso and Ratanov, in preparation). For $\al=1$, by contrast, $(\sT(t),\ep(t))$ is a Markov process.
\end{remark}

\subsection{Special functions}\label{sec:specfun}

We use the confluent hypergeometric \textup(Kummer\textup) function ${}_1F_1$ and the
two-parameter Mittag-Leffler function $E_{\rho,\nu}$,
\begin{equation}
\label{eq:1F1}
{}_1F_1(a,b;z)=\sum_{m\ge0}\frac{(a)_m}{(b)_m}\frac{z^m}{m!},
\qquad
E_{\rho,\nu}(z)=\sum_{m\ge0}\frac{z^m}{\Gamma(\nu+\rho m)},
\end{equation}
where $(a)_m=\Gamma(a+m)/\Gamma(a)$ is the Pochhammer symbol\textup; for the
Mittag-Leffler function and its applications see \citet{saxena}.
Two classical facts about ${}_1F_1$ will be needed: Kummer's transformation
\citep[13.1.27]{AS},
\begin{equation}
\label{eq:kummer}
\rme^{-z}\,{}_1F_1(a,b;z)={}_1F_1(b-a,b;-z),
\end{equation}
and, for $a\ne0,-1,-2,\dots$, the asymptotic behaviour \citep[13.5.1]{AS}
\begin{equation}
\label{eq:1F1asym}
{}_1F_1(a,b;z)\sim\frac{\Gamma(b)}{\Gamma(a)}\,z^{a-b}\rme^{z},
\qquad z\to+\infty.
\end{equation}
The Laplace transform pairs used below, all of the form
$\tau^{1-\nu}/(\tau+b)$ or obtained from it by partial fractions, are tabulated in
\citet[2.1.2.71]{PBM5}.

\section{The fractional Kac--Ornstein--Uhlenbeck process}\label{sec:fKOU}
\setcounter{equation}{0}

\subsection{Well-posedness}

Equation \eqref{eq:fKOU} is a linear Volterra equation which we solve pathwise,
that is, separately for each realization of $\ep$.

\begin{proposition}[Well-posedness and variation of constants]\label{prop:wellposed}
Let $0<\al\le1$\textup, $\beta\ge0$ and $x\in\RR$.
For every realization of $\ep$\textup, equation \eqref{eq:fKOU} has exactly one
solution $X^{(\al)}\in C[0,\infty)$\textup; it is adapted to the filtration
$(\sG_t)$ generated by $\ep$\textup, and it is given by
\begin{equation}\label{eq:X-alt}
X^{(\al)}(t)=x\rme^{-\beta t}+\sF^{(\al)}(t)
-\beta\int_0^t\rme^{-\beta(t-s)}\,\sF^{(\al)}(s)\rmd s .
\end{equation}
\end{proposition}

\begin{proof}
The unknown function enters equation \eqref{eq:fKOU}
only through its primitive, so the integral equation is in fact an ordinary differential equation in disguise.

Put $V(t):=\int_0^tX(s)\rmd s$.
We claim that $X$ solves \eqref{eq:fKOU} if and only if $V$ solves
\begin{equation}\label{eq:Vode}
V'(t)=x-\beta V(t)+\sF^{(\al)}(t),\qquad V(0)=0 .
\end{equation}
Indeed, if $X$ solves \eqref{eq:fKOU} then $V\in C^1$ with $V'=X$, and
\eqref{eq:fKOU} is precisely \eqref{eq:Vode}.
Conversely, if $V$ solves \eqref{eq:Vode}, set $X:=V'$; the right-hand side of
\eqref{eq:Vode} is continuous by Lemma \ref{lem:pathF}, so $X\in C[0,\infty)$,
and $\int_0^tX(s)\rmd s=V(t)-V(0)=V(t)$ turns \eqref{eq:Vode} back into \eqref{eq:fKOU}.

Equation \eqref{eq:Vode} is a linear first-order equation with a continuous
inhomogeneous term; its unique solution is
\[
V(t)=\int_0^t\rme^{-\beta(t-s)}\bigl(x+\sF^{(\al)}(s)\bigr)\rmd s .
\]
By the claim, \eqref{eq:fKOU} therefore has exactly one continuous solution,
$X^{(\al)}=x-\beta V+\sF^{(\al)}$, and since
$\beta\int_0^t\rme^{-\beta(t-s)}x\,\rmd s=x\bigl(1-\rme^{-\beta t}\bigr)$,
this is \eqref{eq:X-alt}.

Finally, $\sF^{(\al)}(t)$ is $\sG_t$-measurable by \eqref{eq:fKOU}, and its
paths are continuous; hence the integral in \eqref{eq:X-alt}, being a limit of
Riemann sums, is $\sG_t$-measurable as well, and $X^{(\al)}$ is adapted.
\end{proof}

\subsection{The convolution representation}

The kernels of this paper belong to the two-parameter family
\begin{equation}\label{def:Knu}
K_\nu(t;b):=\frac{t^{\nu-1}}{\Gamma(\nu)}\,{}_1F_1(1,\nu;-bt)
=\sum_{m\ge0}\frac{(-b)^mt^{\nu+m-1}}{\Gamma(\nu+m)},
\qquad t>0,
\end{equation}
defined for $\nu>0$ and $b\ge0$. We abbreviate $K_\nu(t):=K_\nu(t;\beta)$ and
write the second argument only when it differs from $\beta$.
The second parameter is what allows the moments of Section \ref{sec:moments}
to be written in closed form\textup; there it takes the value $b=2\la$.

\begin{theorem}[Representation]\label{theo:repr}
Let $0<\al\le1$\textup, $\beta\ge0$ and $x\in\RR$.
The solution of \eqref{eq:fKOU} is
\begin{equation}\label{eq:repr}
X^{(\al)}(t)=x\rme^{-\beta t}+\int_0^tK_\al(t-s)\,c_{\ep(s)}\rmd s
=x\rme^{-\beta t}+\bigl(K_\al*c_{\ep(\cdot)}\bigr)(t),
\end{equation}
where the memory kernel $K_\al=K_\al(\cdot\,;\beta)$ is
\begin{equation}\label{def:Kal}
K_\al(t)=\frac{t^{\al-1}}{\Gamma(\al)}\,{}_1F_1(1,\al;-\beta t)
=\frac{t^{\al-1}\rme^{-\beta t}}{\Gamma(\al)}\,{}_1F_1(\al-1,\al;\beta t)
=\sum_{m=0}^\infty\frac{(-\beta)^m t^{\al+m-1}}{\Gamma(\al+m)},
\qquad t>0 .
\end{equation}
Equivalently\textup, $K_\al(t)=t^{\al-1}E_{1,\al}(-\beta t)$.
\end{theorem}

\begin{proof}
Write $k_0(t):=t^{\al-1}/\Gamma(\al)$ and $e_\beta(t):=\rme^{-\beta t}$, so that
$\sF^{(\al)}=k_0*c_{\ep(\cdot)}$ (see \eqref{eq:fKOU}). By \eqref{eq:X-alt},
\[
X^{(\al)}(t)-x\rme^{-\beta t}
=\sF^{(\al)}(t)-\beta\bigl(e_\beta*\sF^{(\al)}\bigr)(t)
=\bigl((k_0-\beta\cdot e_\beta*k_0)*c_{\ep(\cdot)}\bigr)(t),
\]
the rearrangement of the convolutions being justified by Fubini's theorem,
since $k_0$ is locally integrable and $c_{\ep}$ is bounded. It remains to
identify the kernel. Integrating the exponential series term by term,
\[
(e_\beta*k_0)(t)
=\sum_{m\ge0}\frac{(-\beta)^m}{m!}\int_0^t(t-u)^m\frac{u^{\al-1}}{\Gamma(\al)}\rmd u
=\sum_{m\ge0}\frac{(-\beta)^mt^{\al+m}}{\Gamma(\al+m+1)}=K_{1+\al}(t),
\]
and therefore
\[
k_0(t)-\beta\,(e_\beta*k_0)(t)
=\sum_{m\ge0}\frac{(-\beta)^mt^{\al+m-1}}{\Gamma(\al+m)}=K_\al(t).
\]
This is the series in \eqref{def:Kal}. It equals the first expression there by
\eqref{eq:1F1}, the second by Kummer's transformation \eqref{eq:kummer}, and
the Mittag-Leffler form again by \eqref{eq:1F1}.
\end{proof}

\begin{corollary}[The case $\al=1$]\label{cor:alpha1}
For $\al=1$\textup, $K_1(t)=\rme^{-\beta t}$ and \eqref{eq:repr} becomes
\begin{equation}\label{eq:KOUsol36}
X^{(1)}(t)=x\rme^{-\beta t}+\int_0^t\rme^{-\beta(t-s)}c_{\ep(s)}\rmd s,
\end{equation}
the Kac--Ornstein--Uhlenbeck process of \citet{KOU1,KOU2}.
It is the unique solution of the initial value problem
$\rmd X^{(1)}=(c_{\ep(t)}-\beta X^{(1)})\rmd t$\textup, $X^{(1)}(0)=x$.
\end{corollary}

\begin{proof}
Setting $\al=1$ in \eqref{def:Kal} gives
$K_1(t)={}_1F_1(1,1;-\beta t)=\rme^{-\beta t}$, whence \eqref{eq:KOUsol36}
follows. Since $\sF^{(1)}=\sT$ is absolutely continuous with $\sT'=c_{\ep}$,
equation \eqref{eq:fKOU} may be differentiated, which turns it into the stated
initial value problem; uniqueness again follows from Proposition \ref{prop:wellposed}.
\end{proof}

Between switchings $X^{(1)}$ follows the deterministic flow of
$X'=c_j-\beta X$, so that the pair $(X^{(1)},\ep)$ is a piecewise deterministic
Markov process, whose generator is the sum of the generator $\Lambda$ of $\ep$
and the transport operators $(c_j-\beta y)\,\partial_y$ along the two flows; see
\citet[\S2]{KOU1} for the two flows and \citet[\S3]{KOU2} for the generator.

We record for later use some elementary identities of the family
\eqref{def:Knu}.

\begin{lemma}[Elementary identities]\label{lem:kernel-id}
Let $\nu>0$ and $b\ge0$. Then\textup, for $t>0$\textup:
\begin{enumerate}
\item[\textup{(i)}]
\begin{equation}\label{eq:Kint}
\int_0^tK_\nu(u;b)\rmd u=K_{1+\nu}(t;b);
\end{equation}
\item[\textup{(ii)}]
\begin{equation}\label{eq:Kode}
K_\nu(t;b)=\frac{t^{\nu-1}}{\Gamma(\nu)}-b\,K_{1+\nu}(t;b);
\end{equation}
\item[\textup{(iii)}]
\begin{equation}\label{eq:Kvarconst}
K_{1+\nu}(t;b)=\frac{1}{\Gamma(\nu)}\int_0^t\rme^{-b(t-u)}u^{\nu-1}\rmd u>0;
\end{equation}
\item[\textup{(iv)}]
\begin{equation}\label{eq:Kder}
K_\nu'(t;b)=\frac{(\nu-1)\,t^{\nu-2}}{\Gamma(\nu)}-b\,K_\nu(t;b);
\end{equation}
\item[\textup{(v)}] if $b>0$\textup, then
\begin{equation}\label{eq:Kscale}
K_\nu(t;b)=b^{\,1-\nu}\,K_\nu(bt;1),
\end{equation}
so that $b$ is a pure scale parameter and $\nu$ the only shape parameter\textup;
\item[\textup{(vi)}] the Laplace transform of $K_\nu(\cdot\,;b)$ converges
absolutely for $\RRe\,\tau>0$\textup, and
\begin{equation}\label{eq:LTKal}
\widehat{K_\nu(\cdot\,;b)}(\tau):=\int_0^\infty\rme^{-\tau t}K_\nu(t;b)\rmd t
=\frac{\tau^{1-\nu}}{\tau+b},\qquad \RRe\,\tau>0 .
\end{equation}
\end{enumerate}
\end{lemma}

\begin{proof}
\textup{(i)} Integrate the series \eqref{def:Knu} term by term.
\textup{(ii)--(iii)} This is the computation in the proof of Theorem
\ref{theo:repr}, which used nothing about $\al$ and $\beta$ beyond $\nu>0$ and
$b\ge0$.
\textup{(iv)} Differentiate \eqref{eq:Kode} and use \eqref{eq:Kint}.
\textup{(v)} Compare the two series term by term.
\textup{(vi)} By \eqref{eq:Kode} and \eqref{eq:Kvarconst},
$K_\nu(\cdot\,;b)=k-b\,e_b*k$ with $k(t)=t^{\nu-1}/\Gamma(\nu)$ and
$e_b(t)=\rme^{-bt}$. The Laplace integrals of $k$ and of $e_b$ converge
absolutely for $\RRe\,\tau>0$, to $\tau^{-\nu}$ and $1/(\tau+b)$\textup; hence
so does that of $e_b*k$, to their product, and
$\tau^{-\nu}-b\,\tau^{-\nu}/(\tau+b)=\tau^{1-\nu}/(\tau+b)$.
\end{proof}

\begin{lemma}[Paths of $X^{(\al)}$]\label{theo:path}
Let $0<\al\le1$\textup, $\beta\ge0$ and $x\in\RR$. On the event
$\{\ep(0)=i,\ N(t)=n\}$\textup,
\begin{equation}\label{eq:path}
X^{(\al)}(t)=x\rme^{-\beta t}+\Phi_n^{(i)}(t;\tau_1,\dots,\tau_n),
\end{equation}
where\textup, with $\tau_0=0$\textup,
\begin{equation}\label{eq:path2}
\begin{aligned}
\Phi_n^{(i)}(t;\tau_1,\dots,\tau_n)
&:=c_iK_{1+\al}(t)+2c\sum_{k=1}^{n}(-1)^{\ep(\tau_k)}K_{1+\al}(t-\tau_k)\\
&\phantom{:}=\sum_{k=0}^{n}c_{i_k}\Bigl[K_{1+\al}(t-\tau_k)-K_{1+\al}(t-\tau_{k+1}\wedge t)\Bigr].
\end{aligned}
\end{equation}
\end{lemma}

\begin{proof}
By \eqref{eq:repr} it is enough to evaluate the convolution
$(K_\al*c_{\ep(\cdot)})(t)$. On $[\tau_k,\tau_{k+1})$ the velocity equals
$c_{i_k}$, and by \eqref{eq:Kint}
\[
\int_{\tau_k}^{\tau_{k+1}\wedge t}K_\al(t-u)\rmd u
=K_{1+\al}(t-\tau_k)-K_{1+\al}(t-\tau_{k+1}\wedge t).
\]
Summing over $k$ gives the second expression for $\Phi_n^{(i)}$; the first
follows by summation by parts, using
$c_{i_k}-c_{i_{k-1}}=2c(-1)^{\ep(\tau_k)}$ and $K_{1+\al}(0)=0$.
\end{proof}

The representation also gives the exact modulus of continuity of $X^{(\al)}$,
which again separates $\al<1$ sharply from $\al=1$.

\begin{corollary}[Path regularity]\label{cor:holder}
Let $0<\al\le1$\textup, $\beta\ge0$ and $x\in\RR$\textup, and fix $T>0$. There
is a constant $C=C(\al,\beta,T,x,M)$\textup, \emph{not depending on the
realization of $\ep$}\textup, such that
\begin{equation}\label{eq:X-holder}
\bigl|X^{(\al)}(t_2)-X^{(\al)}(t_1)\bigr|\le C\,(t_2-t_1)^{\al},
\qquad 0\le t_1\le t_2\le T.
\end{equation}
For $0<\al<1$ the exponent is sharp\textup: at every switching time $\tau_k$\textup,
\begin{equation}\label{eq:cusp}
X^{(\al)}(\tau_k+\delta)-X^{(\al)}(\tau_k)
=2c\,(-1)^{\ep(\tau_k)}\frac{\delta^\al}{\Gamma(1+\al)}+O(\delta),
\qquad\delta\downarrow0,
\end{equation}
so that $X^{(\al)}$ has an infinite one-sided derivative there. For
$\al=1$\textup, by contrast\textup, $X^{(1)}$ is piecewise $C^\infty$ and its
derivative merely jumps at $\tau_k$\textup, by $2c\,(-1)^{\ep(\tau_k)}$.
\end{corollary}

\begin{proof}
Write $\|f\|_{\infty,T}:=\sup_{[0,T]}|f|$, and estimate the three terms of
\eqref{eq:X-alt}. The first satisfies
$|x||\rme^{-\beta t_2}-\rme^{-\beta t_1}|\le|x|\beta(t_2-t_1)$, and the second
is controlled by \eqref{eq:holder}. For the third, set
\[
g(t):=\int_0^t\rme^{-\beta(t-s)}\sF^{(\al)}(s)\rmd s;
\]
for $0\le t_1\le t_2\le T$,
\begin{align*}
g(t_2)-g(t_1)&=\int_0^{t_1}\left(\rme^{-\beta(t_2-s)}-\rme^{-\beta(t_1-s)}\right)\sF^{(\al)}(s)\rmd s
+\int_{t_1}^{t_2}\rme^{-\beta(t_2-s)}\sF^{(\al)}(s)\rmd s\\
&=\bigl(\rme^{-\beta(t_2-t_1)}-1\bigr)g(t_1)
+\int_{t_1}^{t_2}\rme^{-\beta(t_2-s)}\sF^{(\al)}(s)\rmd s,
\end{align*}
hence
\[
\beta|g(t_2)-g(t_1)|\le\beta(1-\rme^{-\beta t_1})\|\sF^{(\al)}\|_{\infty,T}(t_2-t_1)+\beta\|\sF^{(\al)}\|_{\infty,T}(t_2-t_1)
\le2\beta\|\sF^{(\al)}\|_{\infty,T}(t_2-t_1).
\]
Since $\|\sF^{(\al)}\|_{\infty,T}\le MT^\al/\Gamma(1+\al)$ and
$t_2-t_1\le T^{1-\al}(t_2-t_1)^\al$, we obtain \eqref{eq:X-holder} with
\[
C=\frac{2M}{\Gamma(1+\al)}
+\beta T^{1-\al}\Bigl(|x|+\frac{2MT^\al}{\Gamma(1+\al)}\Bigr).
\]
None of these quantities depends on $\ep$.

For the expansion at $\tau_k$, the single term of \eqref{eq:path2} that
carries $\tau_k$ contributes $2c(-1)^{\ep(\tau_k)}K_{1+\al}(\delta)$, and by
\eqref{def:Knu}
\[
K_{1+\al}(\delta)=\frac{\delta^\al}{\Gamma(1+\al)}+O(\delta^{1+\al}),
\]
while every other term is differentiable at $\tau_k$.
The statement for $\al=1$ follows from Corollary \ref{cor:alpha1}.
\end{proof}

\begin{figure}[ht]
\begin{center}
\includegraphics[width=\textwidth]{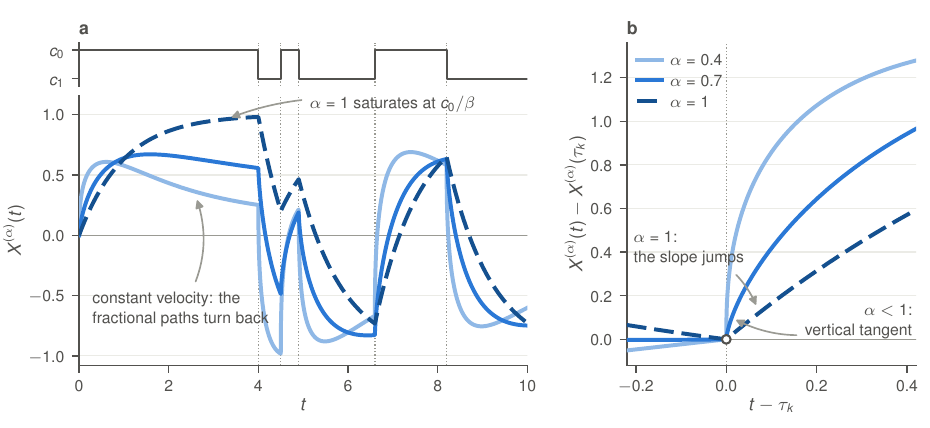}
\end{center}
\caption{Three paths of $X^{(\al)}$ driven by the same velocity path, for
$\al\in\{0.4,\,0.7,\,1\}$ \textup(lighter to darker\textup; the case $\al=1$
is dashed\textup), with $\beta=1$, $c_0=-c_1=1$, $x=0$ and $\ep(0)=0$. The
switching times are $\tau_1,\dots,\tau_5=4.0,\,4.5,\,4.9,\,6.6,\,8.2$
\textup(dotted lines\textup).
\emph{Panel \textup{(a)}}\textup: the velocity $c_{\ep(t)}$ \textup(top\textup)
and the paths, computed from \eqref{eq:path}--\eqref{eq:path2}. Before the
first switching the velocity is constant and $X^{(\al)}(t)=c_0K_{1+\al}(t)$:
for $\al<1$ the path rises only until $t^*$ \textup(here $t^*=0.607$ and
$1.582$\textup) and then turns back towards the origin, by Theorem
\ref{theo:kernel}\,\textup{(v)}, whereas for $\al=1$ it increases towards
$c_0/\beta$.
\emph{Panel \textup{(b)}}\textup: the increments
$X^{(\al)}(t)-X^{(\al)}(\tau_k)$ around the switching time $\tau_k$ with $k=4$,
where the velocity switches from $c_1$ to $c_0$. For $\al<1$ the path leaves
$\tau_k$ with a vertical tangent, growing like $\delta^\al$ by \eqref{eq:cusp};
for $\al=1$ only the slope jumps \textup(Corollary \ref{cor:holder}\textup)}
\label{fig:process}
\end{figure}

\begin{remark}[Cusps]\label{rem:cusp}
For $0<\al<1$ the paths have, by \eqref{eq:cusp}, cusp-type singular points:
each switching time corresponds to a point with a vertical tangent
\textup(Figure \ref{fig:process}\,\textup{(b)}\textup). These
points also explain why estimate \eqref{eq:X-holder}, while uniform with
respect to paths, is of no use for the scaling limit of Section \ref{sec:kac}: there $M=c\to\infty$, so that $C\to\infty$.
\end{remark}

\section{The memory kernel}\label{sec:kernel}
\setcounter{equation}{0}

This section is the technical heart of the paper: every qualitative difference
between $\al<1$ and $\al=1$ can be traced back to Theorem \ref{theo:kernel}.

\subsection{The shape of $K_\al$}
\begin{theorem}[Shape of the kernel]\label{theo:kernel}
Let $\beta>0$.
\begin{enumerate}
\item[\textup{(i)}] If $\al=1$\textup, then $K_1(t)=\rme^{-\beta t}>0$ for all
$t>0$\textup, and $\int_0^\infty K_1(t)\rmd t=1/\beta$.
\end{enumerate}
In the remaining items let $0<\al<1$.
\begin{enumerate}
\item[\textup{(ii)}] \textup{(Behaviour at the ends).}
$K_\al(t)\sim t^{\al-1}/\Gamma(\al)$ as $t\downarrow0$\textup, and\textup, as
$t\to\infty$\textup,
\begin{align}
\label{eq:Kasym}
K_\al(t)&\sim\frac{\al-1}{\beta\,\Gamma(\al)}\,t^{\al-2},\\
\label{eq:Kasym+1}
K_{1+\al}(t)&\sim\frac{t^{\al-1}}{\beta\,\Gamma(\al)}.
\end{align}
\item[\textup{(iii)}] \textup{(One change of sign).} The function
$t\mapsto\rme^{\beta t}K_\al(t)$ is strictly decreasing on $(0,\infty)$\textup,
from $+\infty$ to $-\infty$. Consequently $K_\al$ has a unique zero
$t^*=t^*(\al,\beta)\in(0,\infty)$\textup, with $K_\al>0$ on $(0,t^*)$ and
$K_\al<0$ on $(t^*,\infty)$.
\item[\textup{(iv)}] \textup{(Zero total mass).} $K_\al\in L^1(0,\infty)$ and
\begin{equation}\label{eq:Kzeromass}
\int_0^\infty K_\al(u)\rmd u=0 .
\end{equation}
Moreover\textup,
\begin{equation}\label{eq:Katail}
\int_t^\infty\bigl|K_\al(u)\bigr|\rmd u=K_{1+\al}(t),\qquad t\ge t^*,
\end{equation}
and
\begin{equation}\label{eq:KL1}
\|K_\al\|_{L^1(0,\infty)}=2K_{1+\al}(t^*)
=\frac{2\,(t^*)^{\al-1}}{\beta\,\Gamma(\al)}.
\end{equation}
\item[\textup{(v)}] \textup{(The primitive).} $K_{1+\al}$ is strictly
increasing on $(0,t^*)$ and strictly decreasing on $(t^*,\infty)$\textup; it
attains its global maximum at $t^*$ and tends to $0$ at infinity.
\end{enumerate}
For $\al=1$ we set $t^*:=+\infty$\textup, so that \textup{``}$K_\al>0$ on
$(0,t^*)$\textup{''} remains true\textup; then
$K_2(t)=(1-\rme^{-\beta t})/\beta$ is strictly increasing\textup, with supremum
$1/\beta$ approached only as $t\to\infty$.
\end{theorem}

\begin{proof}
\textup{(i)} is immediate from \eqref{def:Kal}.

\textup{(ii)} The behaviour at the origin is read off the series
\eqref{def:Kal}. At infinity, the second expression in \eqref{def:Kal} and the
asymptotic formula \eqref{eq:1F1asym}, with $a=\al-1$, $b=\al$ and $z=\beta t$,
give
\[
K_\al(t)=\frac{t^{\al-1}\rme^{-\beta t}}{\Gamma(\al)}\,{}_1F_1(\al-1,\al;\beta t)
\sim\frac{t^{\al-1}}{\Gamma(\al-1)}\,(\beta t)^{-1}
=\frac{\al-1}{\beta\,\Gamma(\al)}\,t^{\al-2},
\]
which is \eqref{eq:Kasym}. By \eqref{eq:Kode},
$K_{1+\al}(t)=\beta^{-1}\bigl[t^{\al-1}/\Gamma(\al)-K_\al(t)\bigr]$, where the
second term is of lower order\textup; this is \eqref{eq:Kasym+1}.

\textup{(iii)} By \eqref{eq:Kder},
\[
\frac{\rmd}{\rmd t}\Bigl(\rme^{\beta t}K_\al(t)\Bigr)
=\frac{\al-1}{\Gamma(\al)}\,t^{\al-2}\rme^{\beta t}<0,\qquad t>0 .
\]
By \textup{(ii)}, $\rme^{\beta t}K_\al(t)\to+\infty$ as $t\downarrow0$ and
$\rme^{\beta t}K_\al(t)\to-\infty$ as $t\to\infty$. Hence $\rme^{\beta t}K_\al$
has exactly one zero, which is the unique zero of $K_\al$, and the signs follow.

\textup{(iv)} By \textup{(ii)}, $K_\al$ is integrable at the origin, where it
behaves like $t^{\al-1}$, and at infinity, where it behaves like $t^{\al-2}$.
By \eqref{eq:Kint} and \eqref{eq:Kasym+1},
$\int_0^tK_\al(u)\rmd u=K_{1+\al}(t)\to0$, which is \eqref{eq:Kzeromass}.
Since $K_\al<0$ on $(t^*,\infty)$, it follows that for $t\ge t^*$
\[
\int_t^\infty|K_\al|=-\int_t^\infty K_\al=\int_0^tK_\al=K_{1+\al}(t),
\]
which is \eqref{eq:Katail}. Taking $t=t^*$ and adding
$\int_0^{t^*}|K_\al|=K_{1+\al}(t^*)$ gives the first expression in
\eqref{eq:KL1}\textup; the second is \eqref{eq:Kode} at $t=t^*$, where
$K_\al(t^*)=0$.

\textup{(v)} Since $K_{1+\al}'=K_\al$ by \eqref{eq:Kint}, this follows from
\textup{(iii)} and \eqref{eq:Kasym+1}. The statements for $\al=1$ are
immediate.
\end{proof}

\begin{figure}[ht]
\begin{center}
\includegraphics[width=\textwidth]{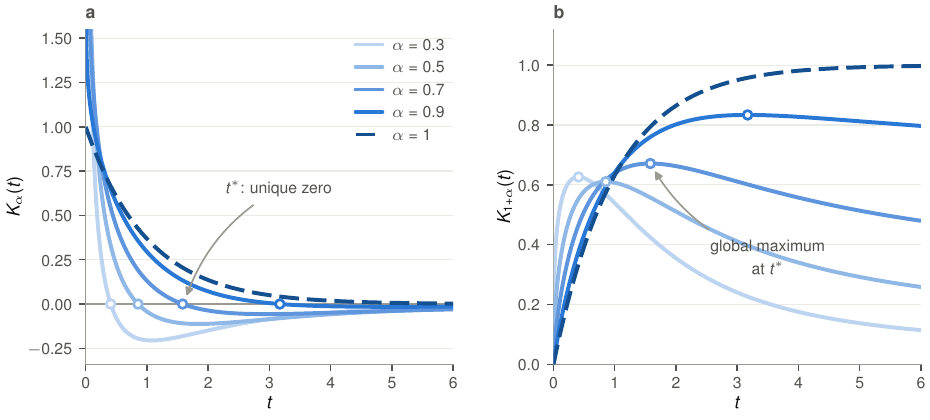}
\end{center}
\caption{The memory kernel and its primitive, for $\beta=1$ and
$\al\in\{0.3,\,0.5,\,0.7,\,0.9,\,1\}$ \textup(lighter to darker\textup); the
case $\al=1$ is dashed.
\emph{Panel \textup{(a)}}\textup: $K_\al$, illustrating Theorem
\ref{theo:kernel}. Each kernel with $\al<1$ descends from $+\infty$, crosses
zero exactly once at $t^*$ \textup(circles\textup), and returns to $0$ from
below; $K_1=\rme^{-\beta t}$ stays positive.
Here $t^*=0.408,\,0.854,\,1.582,\,3.170$; numerically $t^*$ increases with
$\al$ and diverges as $\al\uparrow1$, and by Corollary \ref{cor:scaling} it
scales as $\beta^{-1}$.
\emph{Panel \textup{(b)}}\textup: the primitive $K_{1+\al}$, which by
\eqref{eq:Kint} attains its global maximum exactly at $t^*$ and then decays to
zero, since $\int_0^\infty K_\al=0$ by \eqref{eq:Kzeromass}.
The contrast with $\al=1$, where $K_2\uparrow1/\beta$, is exactly what
Remark \ref{rem:kernel-interp} expresses: the fractional Kac--Ornstein--Uhlenbeck process has zero total memory for
$\al<1$ and positive total memory for $\al=1$
}
\label{fig:kernel}
\end{figure}

\begin{remark}[Zero total memory]\label{rem:kernel-interp}
Identity \eqref{eq:Kzeromass} says that the fractional process $X^{(\al)}$ has \emph{zero total
memory}: a unit impulse of drift at time $0$ has no net effect on $X^{(\al)}$ in
the long run, its immediate positive influence being exactly cancelled by a
long negative tail. Figure \ref{fig:process}\,\textup{(a)} shows the effect on
a path: under a constant velocity the fractional paths turn back towards the
origin.
This is a genuinely fractional phenomenon: for $\al=1$ the total memory is
$1/\beta>0$.

Figure \ref{fig:kernel} shows both kernels for several values of $\al$.
\end{remark}

\begin{corollary}[Scaling of $t^*$]\label{cor:scaling}
Let $0<\al<1$ and write $t^*(\al,\beta)$ for the zero of $K_\al(\cdot\,;\beta)$
given by Theorem \textup{\ref{theo:kernel}}.
Then
\begin{equation}\label{eq:tstar-scale}
t^*(\al,\beta)=\frac{\theta(\al)}{\beta},
\qquad
K_{1+\al}\bigl(t^*;\beta\bigr)=\beta^{-\al}\,K_{1+\al}\bigl(\theta(\al);1\bigr),
\end{equation}
where $\theta(\al):=t^*(\al,1)$ depends on $\al$ alone.
\end{corollary}

\begin{proof}
By \eqref{eq:Kscale}, $K_\al(t;\beta)=0$ if and only if
$K_\al(\beta t;1)=0$, which by uniqueness of the zero means $\beta t=\theta(\al)$;
and \eqref{eq:Kscale} with $\nu=1+\al$ gives the second identity.
\end{proof}

\subsection{The invariant band}

\begin{theorem}[Invariant band]\label{theo:band}
Let $0<\al\le1$\textup, $\beta>0$ and $x\in\RR$\textup, let $t^*$ be as in
Theorem \textup{\ref{theo:kernel}}\textup, and put
\begin{gather}\label{def:Z}
Z(t):=X^{(\al)}(t)-x\rme^{-\beta t}=\int_0^tK_\al(t-u)\,c_{\ep(u)}\rmd u, \\
\label{eq:band}
m_\al(t):=c_0K_{1+\al}(t)-2cK_{1+\al}(t\wedge t^*),
\qquad
M_\al(t):=c_1K_{1+\al}(t)+2cK_{1+\al}(t\wedge t^*).
\end{gather}
\begin{enumerate}
\item[\textup{(i)}] For every $t>0$\textup,
\[
m_\al(t)\le Z(t)\le M_\al(t);
\]
for $t\le t^*$ the two bounds are $c_1K_{1+\al}(t)$ and $c_0K_{1+\al}(t)$.
\item[\textup{(ii)}] The bounds are sharp\textup: $M_\al(t)$ and $m_\al(t)$ are
the supremum and the infimum of $\int_0^tK_\al(t-u)\,c_{\ep(u)}\rmd u$ over all
piecewise constant paths $\ep\colon[0,t]\to\{0,1\}$. For $t>t^*$ the supremum is
attained by the path with exactly one switching\textup, at time
$t-t^*$\textup, from state $1$ to state $0$\textup, and the infimum by its
mirror image\textup; for $t\le t^*$ they are attained by the constant paths.
\item[\textup{(iii)}] Let $0<\al<1$. As $t\to\infty$\textup,
\[
m_\al(t)\ \longrightarrow\ -2cK_{1+\al}(t^*),
\qquad
M_\al(t)\ \longrightarrow\ 2cK_{1+\al}(t^*)=c\,\|K_\al\|_{L^1(0, \infty)},
\]
so that the limiting band is \emph{symmetric} about the origin and depends on
$c_0,c_1$ only through $2c=c_0-c_1$.
If moreover $c_1\le0\le c_0$\textup, then $m_\al$ is nonincreasing and $M_\al$
is nondecreasing on $(t^*,\infty)$\textup, and
\[
Z(t)\in\bigl[-2cK_{1+\al}(t^*),\ 2cK_{1+\al}(t^*)\bigr]
\qquad\text{for every }t>0 ,
\]
the inclusion being into the open interval when $c_1<0<c_0$.
\end{enumerate}
\end{theorem}

\begin{proof}
Representation \eqref{def:Z} follows from \eqref{eq:repr}.

For $t\le t^*$ the kernel $u\mapsto K_\al(t-u)$ is nonnegative on $(0,t)$, with
integral $K_{1+\al}(t)$ by \eqref{eq:Kint}; hence
$c_1K_{1+\al}(t)\le Z(t)\le c_0K_{1+\al}(t)$, with equality holding for the ballistic (switch-free)
paths.
For $t>t^*,$ the kernel is positive exactly on $(t-t^*,t)$, where its integral
is $K_{1+\al}(t^*)$, and negative on $(0,t-t^*)$, where its integral is
$K_{1+\al}(t)-K_{1+\al}(t^*)<0$. The integral is therefore largest when
$c_{\ep}=c_0$ on the first set and $c_{\ep}=c_1$ on the second, where it equals
\[
c_0K_{1+\al}(t^*)+c_1\bigl[K_{1+\al}(t)-K_{1+\al}(t^*)\bigr]
=c_1K_{1+\al}(t)+2cK_{1+\al}(t^*)=M_\al(t);
\]
the corresponding path has exactly one switching, at $t-t^*$. The infimum is
treated symmetrically. This proves \textup{(i)} and \textup{(ii)}.

For \textup{(iii)}, $K_{1+\al}(t)\to0$ by \eqref{eq:Kasym+1}, which gives the
two limits whatever the signs of $c_0,c_1$\textup; by \eqref{eq:KL1} the
limit $2cK_{1+\al}(t^*)$ equals $c\,\|K_\al\|_{L^1}$. By Theorem
\ref{theo:kernel}\,\textup{(v)}, $K_{1+\al}$ is strictly decreasing on
$(t^*,\infty)$, so if $c_1\le0\le c_0$ then $c_0K_{1+\al}(t)\downarrow0$ and
$c_1K_{1+\al}(t)\uparrow0$ there; this is the monotonicity, and with
\textup{(i)} it gives the inclusion for $t>t^*$. For $t\le t^*$ the inclusion
follows from $0\le K_{1+\al}(t)\le K_{1+\al}(t^*)$ and $\max(c_0,-c_1)\le2c$.
\end{proof}

The extremal paths from \textup{(ii)} have probability zero. The fact that $m_\al(t)$ and
$M_\al(t)$ are nonetheless boundary points of the support of the distribution of $Z(t)$,
regardless of the initial state, is proved in Proposition \ref{prop:support}.

\begin{corollary}[The case $\al=1$]\label{cor:attractor}
For $\al=1$ one has $t^*=+\infty$ and $K_2(t)=(1-\rme^{-\beta t})/\beta$\textup,
so that Theorem \textup{\ref{theo:band}\,(i)} reads\textup, with
$\rho_j:=c_j/\beta$\textup,
\[
\rho_1+(x-\rho_1)\rme^{-\beta t}\le X^{(1)}(t)\le\rho_0+(x-\rho_0)\rme^{-\beta t},
\qquad t>0 .
\]
Letting $t\to\infty$ recovers the attractor $[\rho_1,\rho_0]$ of \citet{KOU1}.
\end{corollary}

\begin{figure}[ht]
\begin{center}
\includegraphics[width=\textwidth]{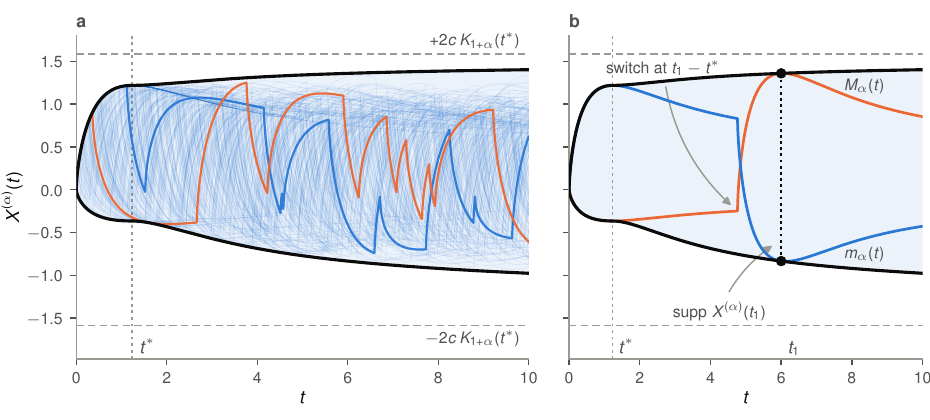}
\end{center}
\caption{The invariant band, for $\al=0.65$, $\beta=1.1$, $c_0=2$, $c_1=-0.6$,
$\la_0=0.8$, $\la_1=1.3$ and $x=0$, so that $t^*=1.2354$,
$K_{1+\al}(t^*)=0.6097$ and $2cK_{1+\al}(t^*)=1.5851$.
In both panels the two heavy curves are $m_\al(t)$ and $M_\al(t)$ of
\eqref{eq:band}, the shaded region between them is the support of the law of
$X^{(\al)}(t)$ \textup(Proposition \ref{prop:support}\textup), and the horizontal
dashed lines are the asymptotic band $\pm2cK_{1+\al}(t^*)$ of Theorem
\ref{theo:band}\textup{(iii)}.
\emph{Panel \textup{(a)}}\textup: $200$ sample paths started at $\ep(0)=0$, two
of them highlighted. Every path leaves the origin \emph{along} the upper
envelope\textup: while $t<\tau_1\wedge t^*$ one has
$X^{(\al)}(t)=c_0K_{1+\al}(t)$, which by \eqref{eq:band} is exactly $M_\al(t)$
for $t\le t^*$. A path that has not switched by $t^*$ detaches from the
envelope at $t^*$, not at $\tau_1$.
\emph{Panel \textup{(b)}}\textup: the two extremal paths associated with a
single instant $t_1=6$. Each switches exactly once, at $t_1-t^*$, and touches
the corresponding envelope \emph{at $t_1$ only}\textup; by Theorem
\ref{theo:band}\textup{(ii)} they start in \emph{opposite} states, the upper
one at $\ep(0)=1$ and the lower one at $\ep(0)=0$. The vertical dotted segment
they span is the support of $X^{(\al)}(t_1)$ \textup(Proposition
\ref{prop:support}\textup).
Note that the envelopes approach the asymptotic band slowly, like
$t^{\al-1}$\textup: the width of the gap at $t=10$ is still $11.5\%$ of
$2cK_{1+\al}(t^*)$ above and $38.4\%$ below. The band is symmetric in the limit
but the approach is not, because $|c_0|\ne|c_1|$}
\label{fig:band}
\end{figure}

\begin{remark}
\label{rem:band-contrast}
For $\al=1$ the limiting band $[\rho_1,\rho_0]$ is asymmetric, unless
$c_0+c_1=0$, and absorbing: by Corollary \ref{cor:attractor} and the Markov
property, once the process enters it, it never leaves.
For $\al<1$ the endpoints of $[m_\al(t),M_\al(t)]$ converge to the
\emph{symmetric} limits $\mp2cK_{1+\al}(t^*)$, which do not depend on the signs
of $c_0$ and $c_1$; see Figure \ref{fig:band}.
\end{remark}

\begin{proposition}
\label{prop:support}
Let $0<\al<1$\textup, $\beta>0$\textup, $x\in\RR$\textup, $t>0$ and
$i\in\{0,1\}$. Under $\PP_i,$ the support of the law of $Z(t)$ is
$[m_\al(t),M_\al(t)]$\textup; that is\textup,
\[
\supp X^{(\al)}(t)=x\rme^{-\beta t}+\bigl[m_\al(t),M_\al(t)\bigr].
\]
In particular\textup, the support does not depend on the initial state.
\end{proposition}

\begin{proof}
The inclusion $\supp Z(t)\subseteq[m_\al(t),M_\al(t)]$ follows from Theorem
\ref{theo:band}\,\textup{(i)}, the interval being closed.

To prove the reverse inclusion, consider the function $\Phi_2^{(i)}(t;\cdot,\cdot)$ defined as in
\eqref{eq:path2}. As a function of $(\tau_1,\tau_2),$ it is an
alternating sum of shifts of $K_{1+\al}$,
\begin{equation}\label{eq:Phi2}
\Phi_2^{(i)}(t;\tau_1,\tau_2)
=c_iK_{1+\al}(t)+2c(-1)^{i+1}\bigl[K_{1+\al}(t-\tau_1)-K_{1+\al}(t-\tau_2)\bigr],
\end{equation}
and is therefore continuous on the closure $\overline{\Delta_2(t)}
=\{(\tau_1,\tau_2):0\le\tau_1\le\tau_2\le t\}$ of the simplex
\eqref{def:simplex}, which is compact and connected. Hence,
\[
R^{(i)}:=\Phi_2^{(i)}\bigl(t;\overline{\Delta_2(t)}\bigr)
\]
is a compact interval, and $R^{(i)}\subseteq[m_\al(t),M_\al(t)]$ by Theorem
\ref{theo:band}\,\textup{(ii)}.

\emph{Both endpoints belong to $R^{(i)}$.}
Put $s^*:=(t-t^*)^+$, so that $K_{1+\al}(t-s^*)=K_{1+\al}(t\wedge t^*)$.
Evaluating \eqref{eq:Phi2} at the four corner configurations yields
\[
\Phi_2^{(0)}(t;0,s^*)=\Phi_2^{(1)}(t;s^*,t)=M_\al(t),
\qquad
\Phi_2^{(0)}(t;s^*,t)=\Phi_2^{(1)}(t;0,s^*)=m_\al(t).
\]
Indeed, for $i=0$ and $(\tau_1,\tau_2)=(0,s^*)$ formula \eqref{eq:Phi2}
gives $c_0K_{1+\al}(t)-2cK_{1+\al}(t)+2cK_{1+\al}(t\wedge t^*)
=c_1K_{1+\al}(t)+2cK_{1+\al}(t\wedge t^*)=M_\al(t)$, and the other three are
analogous computations. \textup(The configurations with $\tau_2=t$ are the
one-switching paths of Theorem \ref{theo:band}\,\textup{(ii)}\textup; those
with $\tau_1=0$ are the same paths started in the opposite state, which is why
both initial states reach both endpoints.\textup)
A compact interval contained in $[m_\al(t),M_\al(t)]$ and containing both of
its endpoints equals it, so
\begin{equation}\label{eq:R2}
R^{(i)}=[m_\al(t),M_\al(t)],\qquad i=0,1 .
\end{equation}

\emph{Every point of $R^{(i)}$ is in the support.}
Let $y=\Phi_2^{(i)}(t;\tau_1,\tau_2)$ with
$(\tau_1,\tau_2)\in\overline{\Delta_2(t)}$, and let $\delta>0$. By continuity
there is an open $O\subseteq\RR^2$ containing $(\tau_1,\tau_2)$ with
$|\Phi_2^{(i)}(t;\cdot)-y|<\delta$ on $O\cap\overline{\Delta_2(t)}$. Since
$\Delta_2(t)$ is open and dense in its closure, the set $O\cap\Delta_2(t)$ is
open and nonempty, hence of positive Lebesgue measure; and by Lemma
\ref{lem:switch}\,\textup{(i)} the density $L_2^{(i)}$ is strictly positive
there. Therefore, by Lemma \ref{theo:path},
\[
\PP_i\bigl\{|Z(t)-y|<\delta\bigr\}
\ \ge\ \int_{O\cap\Delta_2(t)}L_2^{(i)}(\tau_1,\tau_2;t)\,
\rmd\tau_1\rmd\tau_2\ >\ 0 .
\]
As $\delta>0$ was arbitrary, $y\in\supp Z(t)$. With \eqref{eq:R2} this gives
$[m_\al(t),M_\al(t)]\subseteq\supp Z(t)$ and completes the proof.
\end{proof}

\section{Moments and second-order structure}\label{sec:moments}
\setcounter{equation}{0}

By the representation \eqref{eq:repr}, every moment of $X^{(\al)}$ is an
integral of the kernel $K_\al$ against the moments of the velocity, which
Lemma \ref{lem:cov} provides. We first compute the conditional means and
variances, then the autocovariance, and finally pass to the stationary regime,
where the zero total mass of $K_\al$ reappears as antipersistence.

We begin with a general convolution lemma.

\begin{lemma}[Moments of a convolution]\label{lem:conv}
Let $K\in L^1(0,T)$ for every $T>0$\textup, and let
$Y(t):=(K*c_{\ep(\cdot)})(t)$\textup, $t\ge0$. Put
\begin{equation}\label{def:J}
\begin{aligned}
J_0(t)&=\int_0^tK(u)\rmd u, \qquad
J_1(t)=\int_0^tK(t-u)\rme^{-2\la u}\rmd u,\\
J_2(t)&=\iint_{[0,t]^2}K(t-u_1)K(t-u_2)\rme^{-2\la|u_1-u_2|}\rmd u_1\rmd u_2, \\
J_3(t)&=\iint_{[0,t]^2}K(t-u_1)K(t-u_2)\rme^{-2\la(u_1\vee u_2)}\rmd u_1\rmd u_2.
\end{aligned}
\end{equation}
Then\textup, for every $t\ge0$\textup,
\begin{equation}\label{eq:MV-general}
\begin{aligned}
\EE_0[Y(t)]&=\kappa J_0(t)+2cp_1^*J_1(t), \qquad
\EE_1[Y(t)]=\kappa J_0(t)-2cp_0^*J_1(t),\\
\Var_0[Y(t)]&=4c^2p_1^*\bigl[(p_1^*-p_0^*)J_3(t)+p_0^*J_2(t)-p_1^*J_1(t)^2\bigr],\\
\Var_1[Y(t)]&=4c^2p_0^*\bigl[(p_0^*-p_1^*)J_3(t)+p_1^*J_2(t)-p_0^*J_1(t)^2\bigr].
\end{aligned}
\end{equation}
\end{lemma}

\begin{proof}
The velocity is bounded and $K$ is integrable on $(0,t)$, so Fubini's theorem
gives $\EE_i[Y(t)]=\int_0^tK(t-s)\,\EE_i[c_{\ep(s)}]\rmd s$ and, for
$t_1,t_2\ge0$,
\begin{equation}\label{eq:cov-conv}
\Cov_i\bigl(Y(t_1),Y(t_2)\bigr)
=\int_{[0, t_1]\times[0, t_2]}
K(t_1-u_1)K(t_2-u_2)
\Cov_i\bigl(c_{\ep(u_1)},c_{\ep(u_2)}\bigr)\rmd u_1\rmd u_2.
\end{equation}
Inserting \eqref{eq:mean01} into the first identity gives the means. Taking
$t_1=t_2=t$ in \eqref{eq:cov-conv} and inserting
\eqref{eq:ACov0}--\eqref{eq:ACov1} gives the variances\textup: the three
exponentials produce $J_3$, $J_2$ and $J_1^2$, in this order.
\end{proof}

\subsection{Means and variances}

Throughout this section, we apply
\eqref{def:J} to $K=K_\al$ and drop the kernel from the notation, writing
simply $J_0,\dots,J_3$\textup; by \eqref{eq:Kint}, $J_0=K_{1+\al}$.

\begin{theorem}[Means and variances]\label{theo:moments}
Let $0<\al\le1$\textup, $\beta\ge0$ and $x\in\RR$. For every $t\ge0$\textup,
\begin{align}
\label{eq:E0X}
\EE_0[X^{(\al)}(t)]&=x\rme^{-\beta t}+\kappa K_{1+\al}(t)+2cp_1^*J_1(t),\\
\label{eq:E1X}
\EE_1[X^{(\al)}(t)]&=x\rme^{-\beta t}+\kappa K_{1+\al}(t)-2cp_0^*J_1(t),\\
\label{eq:V0X}
\Var_0[X^{(\al)}(t)]&=4c^2p_1^*\bigl[(p_1^*-p_0^*)J_3(t)+p_0^*J_2(t)-p_1^*J_1(t)^2\bigr],\\
\label{eq:V1X}
\Var_1[X^{(\al)}(t)]&=4c^2p_0^*\bigl[(p_0^*-p_1^*)J_3(t)+p_1^*J_2(t)-p_0^*J_1(t)^2\bigr],
\end{align}
where $J_1,J_2,J_3$ are the integrals \eqref{def:J} with $K=K_\al$.
Moreover\textup, for $t>0$\textup,
\begin{equation}\label{eq:J1}
J_1(t)
=\begin{cases}
  \dfrac{K_\al(t)-K_\al(t;2\la)}{2\la-\beta},    & \text{if $2\la\ne\beta$}, \\[3mm]
\dfrac{t^\al}{\Gamma(1+\al)}\,{}_1F_1(2,1+\al;-\beta t),
   & \text{if $2\la=\beta$},
\end{cases}
\end{equation}
and
\begin{align}\label{eq:J2}
J_2(t)&=2\int_0^tK_\al(u) J_1(u)\rmd u,\\
\label{eq:J3}
J_3(t)&=2\int_0^t\rme^{-2\la(t-u)}K_\al(u) J_1(u)\rmd u.
\end{align}
\end{theorem}

In particular, each conditional mean is a combination of three Kummer
functions: for $2\la\ne\beta$,
\begin{equation}\label{eq:mean-closed}
\EE_0\bigl[X^{(\al)}(t)\bigr]=x\rme^{-\beta t}+\kappa K_{1+\al}(t)
+2cp_1^*\,\frac{K_\al(t)-K_\al(t;2\la)}{2\la-\beta},
\end{equation}
and similarly for $\EE_1$, with $-2cp_0^*$ in place of $2cp_1^*$.

\begin{proof}
By \eqref{eq:repr}, $X^{(\al)}(t)-x\rme^{-\beta t}=(K_\al*c_{\ep(\cdot)})(t)$, so
the four moment formulae are Lemma \ref{lem:conv} with $K=K_\al$, for which
$J_0=K_{1+\al}$ by \eqref{eq:Kint}\textup; the deterministic term
$x\rme^{-\beta t}$ enters only the means.

\emph{The integral $J_1$.}
By definition $J_1=K_\al*\rme^{-2\la\,\cdot}$, so by \eqref{eq:LTKal}
\[
\widehat{J_1}(\tau)=\frac{\tau^{1-\al}}{(\tau+\beta)(\tau+2\la)},
\qquad \RRe\,\tau>0 .
\]
If $2\la\ne\beta$, then
\[\frac{1}{(\tau+\beta)(\tau+2\la)}
=\frac{1}{2\la-\beta}\left[\frac{1}{\tau+\beta}-\frac{1}{\tau+2\la}\right],\]
so that $\widehat{J_1}$ is the corresponding combination of
$\widehat{K_\al}$ and $\widehat{K_\al(\cdot\,;2\la)}$, and the first case of
\eqref{eq:J1} follows by uniqueness of the Laplace transform.

For the second case of \eqref{eq:J1}, put $j(t;b):=\int_0^tK_\al(t-u)\rme^{-bu}\rmd u$ for $b>0$,
so that $J_1(t)=j(t;2\la)$. 

For fixed $t>0$ the function $b\mapsto j(t;b)$ is continuous by dominated
convergence, and for $b\ne\beta$ the first case gives
$j(t;b)=\bigl(K_\al(t)-K_\al(t;b)\bigr)/(b-\beta)$. Since the series
\eqref{def:Knu} is an entire function of $b$, letting $b\to\beta$ and
differentiating \eqref{def:Knu} term by term gives
\[
j(t;\beta)=-\partial_bK_\al(t;b)\big|_{b=\beta}
=t^\al\sum_{k\ge0}\frac{(k+1)(-\beta t)^k}{\Gamma(1+\al+k)}
=\frac{t^\al}{\Gamma(1+\al)}\,{}_1F_1(2,1+\al;-\beta t),
\]
which corresponds to the second case of \eqref{eq:J1}.

\emph{The integrals $J_2$ and $J_3$.}
Substituting $u_i\to t-u_i$ in \eqref{def:J} and using the symmetry of the
integrand,
\[
J_2(t)=\iint_{[0,t]^2}K_\al(u_1)K_\al(u_2)\rme^{-2\la|u_1-u_2|}\rmd u_1\rmd u_2
=2\int_0^tK_\al(u_1)\rme^{-2\la u_1}
\Bigl[\int_0^{u_1}K_\al(u_2)\rme^{2\la u_2}\rmd u_2\Bigr]\rmd u_1.
\]
Writing the convolution the other way round,
\begin{equation}\label{eq:I1ode}
J_1(u)=\rme^{-2\la u}\int_0^u\rme^{2\la v}K_\al(v)\rmd v,
\end{equation}
so the bracket equals $\rme^{2\la u_1}J_1(u_1)$, the two exponentials cancel,
and \eqref{eq:J2} follows.
The same substitution turns $\rme^{-2\la(u_1\vee u_2)}$ into
$\rme^{-2\la t}\rme^{2\la(u_1\wedge u_2)}$, and the same computation gives
\eqref{eq:J3}.
\end{proof}

\begin{lemma}[Long-time behaviour of $J_1, J_2, J_3$]\label{lem:Irates}
Let $0<\al<1$ and $\beta>0$. Then $K_\al\cdot J_1\in L^1(0, \infty)$ and\textup, as
$t\to\infty,$
\begin{equation}\label{eq:J1as}
J_1(t)\sim\frac{\al-1}{2\la\beta\Gamma(\al)}\,t^{\al-2},
\end{equation}
\begin{equation}\label{eq:J23as}
J_2(t)\longrightarrow2\int_0^\infty K_\al(u) J_1(u)\rmd u,
\qquad
J_3(t)\longrightarrow0.
\end{equation}
\end{lemma}
\begin{proof}
We use the following elementary fact. Let $\mu>0$, let $g$ be integrable
on $(0,T)$ for every $T>0$, and suppose that $g(t)\sim Ct^{-\rho}$ as
$t\to\infty$, with $C\ne0$ and $\rho>0$. Then
\begin{equation}\label{eq:exp-conv}
\int_0^t\rme^{-\mu(t-u)}g(u)\rmd u\sim\frac{g(t)}{\mu},\qquad t\to\infty.
\end{equation}
Indeed, the part of the integral over $(0,t/2)$ is at most
$\rme^{-\mu t/2}\int_0^{t/2}|g(u)|\rmd u$, which is $O(t\,\rme^{-\mu t/2})$
because $g$ is locally integrable and bounded near infinity. In the part over
$(t/2,t)$ put $u=t-w$\textup: it equals
\[
g(t)\int_0^{t/2}\rme^{-\mu w}\,\frac{g(t-w)}{g(t)}\rmd w ,
\]
where the ratio tends to $1$ for each fixed $w$ and is bounded by $2^{\rho+1}$
for all $w\le t/2$ once $t$ is large\textup; by dominated convergence the
integral tends to $1/\mu$.

\emph{The rate of $J_1$.}
By \eqref{def:J}, $J_1$ is the left-hand side of \eqref{eq:exp-conv} with
$g=K_\al$ and $\mu=2\la$, and
$K_\al(t)\sim(\al-1)t^{\al-2}/(\beta\Gamma(\al))$ by \eqref{eq:Kasym}. Hence
$J_1(t)\sim K_\al(t)/(2\la)$, which is \eqref{eq:J1as}.

\emph{Integrability and the limit of $J_2$.}
As $u\downarrow0$, $K_\al(u)\sim u^{\al-1}/\Gamma(\al)$ by Theorem
\ref{theo:kernel}\,\textup{(ii)}, so
$|J_1(u)|\le\int_0^u|K_\al(v)|\rmd v=O(u^\al)$ and
$K_\al(u)J_1(u)=O(u^{2\al-1})$\textup; as $u\to\infty$,
$K_\al(u)J_1(u)=O(u^{2\al-4})$ by \eqref{eq:Kasym} and \eqref{eq:J1as}.
Both exponents are integrable for $0<\al<1$, so
$K_\al\cdot J_1\in L^1(0,\infty)$, and \eqref{eq:J2} gives the limit of $J_2$.

\emph{The limit of $J_3$.}
Splitting the integral \eqref{eq:J3} at $t/2$,
\[
|J_3(t)|\le2\rme^{-\la t}\,\|K_\al J_1\|_{L^1(0,\infty)}
+\frac1\la\sup_{u\ge t/2}|K_\al(u)J_1(u)|\longrightarrow0 ,
\]
because $K_\al\cdot J_1$ is integrable and tends to $0$ at infinity.
\end{proof}

\begin{corollary}[Long-time behaviour of the moments]\label{cor:moments-limit}
Let $0<\al<1$\textup, $\beta>0$ and $x\in\RR$. Then for both $i\in\{0,1\},$
\begin{equation}\label{eq:mean-limit}
\EE_i\bigl[X^{(\al)}(t)\bigr]
=\frac{\kappa}{\beta\Gamma(\al)}\,t^{\al-1}+O\bigl(t^{\al-2}\bigr)
\longrightarrow0,
\qquad t\to\infty,
\end{equation}
and
\begin{equation}\label{eq:var-limit}
\Var_i\bigl[X^{(\al)}(t)\bigr]\longrightarrow
2\sigma^2\int_0^\infty K_\al(u) J_1(u)\rmd u.
\end{equation}

If $\kappa=0$ \textup(in particular\textup, in the symmetric case\textup{),}
then the leading term of \eqref{eq:mean-limit} vanishes\textup, and
\[\EE_0[X^{(\al)}(t)]\sim2cp_1^*\,\frac{\al-1}{2\la\beta\Gamma(\al)}t^{\al-2},\qquad
\EE_1[X^{(\al)}(t)]\sim-2cp_0^*\,\frac{\al-1}{2\la\beta\Gamma(\al)}t^{\al-2}.\]
\end{corollary}

\begin{proof}
The term $x\rme^{-\beta t}$ in \eqref{eq:E0X}--\eqref{eq:E1X} is exponentially
small; moreover
$K_{1+\al}(t)=\bigl(t^{\al-1}/\Gamma(\al)-K_\al(t)\bigr)/\beta
=t^{\al-1}/(\beta\Gamma(\al))+O(t^{\al-2})$
by \eqref{eq:Kode} and \eqref{eq:Kasym}, and
$J_1(t)=O(t^{\al-2})$ by Lemma \ref{lem:Irates}. This gives
\eqref{eq:mean-limit}.

 If $\kappa=0$ the first term is absent, and
\eqref{eq:J1as} gives the last two relations.

In \eqref{eq:V0X}--\eqref{eq:V1X} both $J_3(t)$ and $J_1(t)^2$ tend to $0$ by
Lemma \ref{lem:Irates}, so only the middle term survives, and
$4c^2p_1^*p_0^*=\sigma^2$ turns it into
$2\sigma^2\int_0^\infty K_\al(u) J_1(u)\rmd u$.

Both limits are independent of the initial state $i$, as expected.
\end{proof}

The formulae \eqref{eq:E0X}--\eqref{eq:E1X} have a qualitative consequence
that has no counterpart at $\al=1$\textup; it rests on the sign of $J_1$.

\begin{lemma}[The sign of $J_1$]\label{lem:I1sign}
Let $0<\al<1$ and $\beta>0$. Then $J_1$ has a unique zero
$t_0=t_0(\al,\beta,\la)$ in $(0,\infty)$\textup, it lies in $(t^*,\infty)$\textup,
and
\[
J_1>0\ \text{ on }(0,t_0),\qquad J_1<0\ \text{ on }(t_0,\infty).
\]
\end{lemma}

\begin{proof}
By \eqref{eq:I1ode}, $\bigl(\rme^{2\la t}J_1(t)\bigr)'=\rme^{2\la t}K_\al(t)$,
which is positive on $(0,t^*)$ and negative on $(t^*,\infty)$ by Theorem
\ref{theo:kernel}\,\textup{(iii)}.
So $\rme^{2\la t}J_1(t)$ increases strictly from $0$ on $(0,t^*)$ and decreases
strictly afterwards, and it is eventually negative, because $J_1(t)<0$ for
large $t$ by \eqref{eq:J1as}.
It has therefore exactly one zero, which lies in $(t^*,\infty)$ and is the
unique zero of $J_1$.
\end{proof}

\begin{corollary}[Crossing of the conditional means]\label{cor:means-cross}
In the setting of Theorem \textup{\ref{theo:moments}}, for every $t\ge0$,
\begin{equation}\label{eq:means-diff}
\EE_0\bigl[X^{(\al)}(t)\bigr]-\EE_1\bigl[X^{(\al)}(t)\bigr]=2cJ_1(t).
\end{equation}
If $0<\al<1$ and $\beta>0$, this difference is positive on $(0,t_0)$ and
negative on $(t_0,\infty)$, where $t_0>t^*$ is the unique zero of $J_1$\textup;
if $\al=1$, it is positive for every $t>0$.
\end{corollary}

\begin{proof}
Subtracting \eqref{eq:E1X} from \eqref{eq:E0X} and using $p_0^*+p_1^*=1$ gives
\eqref{eq:means-diff}. For $0<\al<1$ and $\beta>0$ the sign of $J_1$ is given
by Lemma \ref{lem:I1sign}. For $\al=1$,
$J_1(t)=\int_0^t\rme^{-\beta(t-u)}\rme^{-2\la u}\rmd u>0$ by \eqref{def:J}.
\end{proof}

\subsection{Autocovariance}

The two-time covariance has the same three-term shape as the variance, and
reduces to single integrals in the same way.
Besides $J_1$ and $K_{1+\al}$, one new function is needed: the
\emph{incomplete transform}
\begin{equation}\label{def:Theta}
\Theta(s):=\int_0^\infty\rme^{-2\la u}K_\al(s+u)\rmd u,
\qquad s\ge0 .
\end{equation}
The integral converges absolutely, because $K_\al$ is integrable at the origin
and bounded away from it\textup; by \eqref{eq:LTKal},
$\Theta(0)=\widehat{K_\al}(2\la)=(2\la)^{1-\al}/(2\la+\beta)$\textup; and
$\Theta(s)\to0$ as $s\to\infty$, except when $\al=1$ and $\beta=0$.

\begin{theorem}[Autocovariance]\label{theo:acov}
Let $0<\al\le1$\textup, $\beta\ge0$ and $x\in\RR$\textup, let $0\le t_1\le t_2$
and put $h:=t_2-t_1$. Then
\begin{equation}\label{eq:acov}
\Cov_0\bigl(X^{(\al)}(t_1),X^{(\al)}(t_2)\bigr)
=4c^2p_1^*\Bigl[(p_1^*-p_0^*)\,\sI_3(t_1,t_2)+p_0^*\,\sI_2(t_1,t_2)
-p_1^*\,J_1(t_1)J_1(t_2)\Bigr],
\end{equation}
and $\Cov_1$ is obtained by interchanging $p_0^*$ and $p_1^*$\textup; here
\begin{align}
\label{eq:acov2}
\sI_2(t_1,t_2)&=\int_0^{t_1}K_\al(u)\bigl[J_1(h+u)+\Theta(h+u)\bigr]\rmd u
-\Theta(t_2)\,J_1(t_1),\\
\label{eq:acov3}
\sI_3(t_1,t_2)&=K_{1+\al}(t_2)\,J_1(t_1)
+\int_0^{t_1}\rme^{-2\la(t_1-u)}K_\al(u)
\bigl[J_1(h+u)-K_{1+\al}(h+u)\bigr]\rmd u .
\end{align}
For $t_1=t_2=t$\textup, $\sI_2(t,t)=J_2(t)$ and $\sI_3(t,t)=J_3(t)$\textup, and
\eqref{eq:acov} reduces to \eqref{eq:V0X}.
\end{theorem}

\begin{proof}
The deterministic term $x\rme^{-\beta t}$ does not affect covariances, so by
\eqref{eq:repr} and \eqref{eq:cov-conv} with $K=K_\al$,
\[
\Cov_i\bigl(X^{(\al)}(t_1),\;X^{(\al)}(t_2)\bigr)
=\int_{0}^{t_1}\left(\int_0^{t_2}K_\al(t_1-s_1)K_\al(t_2-s_2)
\Cov_i\bigl(c_{\ep(s_1)},c_{\ep(s_2)}\bigr)\rmd s_2\right)\rmd s_1 ,
\]
and \eqref{eq:ACov0} splits this into three pieces.
The piece carrying $\rme^{-2\la(s_1+s_2)}$ factorizes into
$J_1(t_1)J_1(t_2)$, which is the last term of \eqref{eq:acov}. Denote by
$\sI_2$ and $\sI_3$ the pieces carrying $\rme^{-2\la|s_1-s_2|}$ and
$\rme^{-2\la(s_1\vee s_2)}$\textup.

For $t_1=t_2=t$ they are $J_2(t)$ and
$J_3(t)$ by \eqref{def:J}, which is the last assertion. It remains to show that
they equal the right-hand sides of \eqref{eq:acov2} and \eqref{eq:acov3}.

Substitute $u_i=t_i-s_i$, so that $s_1-s_2=-(h+u_1-u_2)$ and
$s_1\vee s_2=t_2-\min(h+u_1,u_2)$.

Consider $\sI_2$ and split the inner integral at $u_2=h+u_1$.
On $\{u_2<h+u_1\}$, the substitution $w=h+u_1-u_2$ gives inner integral
$\int_0^{h+u_1}\rme^{-2\la w}K_\al(h+u_1-w)\rmd w=J_1(h+u_1)$.
On $\{u_2>h+u_1\}$, the substitution $w=u_2-(h+u_1)$ (whose range is
$(0,t_1-u_1)$, because $t_2-h=t_1$) gives
\[\int_0^{t_1-u_1}\rme^{-2\la w}K_\al(h+u_1+w)\rmd w
=\Theta(h+u_1)-\rme^{-2\la(t_1-u_1)}\Theta(t_2).\]
The last step follows by \eqref{def:Theta} and $h+u_1+(t_1-u_1)=t_2$.
Integrating against $K_\al(u_1)$ and using
$\int_0^{t_1}\rme^{-2\la(t_1-u)}K_\al(u)\rmd u=J_1(t_1)$ gives
\eqref{eq:acov2}.

For $\sI_3,$ the same splitting is applied, with $\rme^{2\la\min(h+u_1,u_2)}$ instead
 of $\rme^{-2\la|h+u_1-u_2|}$.
On $\{u_2>h+u_1\}$ the minimum is $h+u_1,$ and the inner integral is
$\int_{h+u_1}^{t_2}K_\al(s)\rmd s=K_{1+\al}(t_2)-K_{1+\al}(h+u_1)$ by \eqref{eq:Kint}.
On $\{u_2<h+u_1\}$ it is $\rme^{2\la(h+u_1)}J_1(h+u_1)$ as before.
Collecting the prefactor $\rme^{-2\la t_2}\rme^{2\la(h+u_1)}
=\rme^{-2\la(t_1-u_1)}$ gives \eqref{eq:acov3}.
\end{proof}

\subsection{The stationary regime}\label{sec:stationary}

For $0<\al<1$ and $\beta>0$ the kernel $K_\al$ is integrable on the whole
half-line, so the process can be started in the infinite past.

\begin{proposition}[Stationary version]\label{prop:stationary}
Let $0<\al<1,\;\beta>0$, and let $\{\ep(t)\}_{t\in\RR}$ be a two-sided
stationary version of Lemma \textup{\ref{lem:switch}\,(iv)}. Set
\begin{equation}\label{def:Xstat}
\widetilde X^{(\al)}(t):=\int_{-\infty}^tK_\al(t-s)\,c_{\ep(s)}\rmd s
=\int_0^\infty K_\al(u)\,c_{\ep(t-u)}\rmd u,
\qquad t\in\RR .
\end{equation}
The integral converges absolutely on every path\textup. Moreover\textup,
\begin{enumerate}
\item[\textup{(i)}] $\widetilde X^{(\al)}$ is strictly stationary\textup,
\begin{equation}\label{eq:Xstat-mean}
\EE\bigl[\widetilde X^{(\al)}(t)\bigr]=\kappa\int_0^\infty K_\al(u)\rmd u=0 ,
\end{equation}
and\textup, regardless of the signs of $c_0$ and $c_1$\textup,
\begin{equation}\label{eq:Xstat-band}
\bigl|\widetilde X^{(\al)}(t)\bigr|\le2cK_{1+\al}(t^*)
\qquad\text{for every }t\in\RR .
\end{equation}
\item[\textup{(ii)}] Let $X^{(\al)}$ be the solution of \eqref{eq:fKOU} with
$X^{(\al)}(0)=x$ driven by the restriction of the same $\ep$ to $[0,\infty)$.
Then\textup, on every path and for every $t\ge t^*$\textup,
\begin{equation}\label{eq:coupling}
\bigl|X^{(\al)}(t)-\widetilde X^{(\al)}(t)\bigr|
\le|x|\rme^{-\beta t}+M\,K_{1+\al}(t)=O\bigl(t^{\al-1}\bigr),
\qquad t\to\infty .
\end{equation}
\item[\textup{(iii)}] Let $x\in\RR$ and $i\in\{0,1\}$. For every $n\ge1$ and
$0\le h_1<\dots<h_n$\textup, under $\PP_i$\textup,
\begin{equation}\label{eq:limlaw}
\bigl(X^{(\al)}(t+h_1),\dots,X^{(\al)}(t+h_n)\bigr)\Longrightarrow
\bigl(\widetilde X^{(\al)}(h_1),\dots,\widetilde X^{(\al)}(h_n)\bigr),
\qquad t\to\infty,
\end{equation}
where $\Longrightarrow$ denotes weak convergence in $\RR^n$.

All mixed moments also converge.
\end{enumerate}
\end{proposition}

\begin{proof}
By Theorem \ref{theo:kernel}\,\textup{(iv)}, $K_\al\in L^1(0,\infty)$\textup; as
$|c_\ep|\le M$, the integral \eqref{def:Xstat} converges absolutely on every
path.

\textup{(i)}
Let $D$ be the space of $\{0,1\}$-valued c\`adl\`ag paths on $\RR$, with the
$\sigma$-field generated by the coordinate maps, let
$(\theta_tw)(\cdot):=w(\cdot+t)$ be the shift, and put
$\Phi(w):=\int_0^\infty K_\al(u)\,c_{w(-u)}\rmd u, \;w\in D$.

The functional $\Phi$ is measurable\textup: it is the uniform limit, as
$a\downarrow0$ and $R\to\infty$, of the truncated integrals over $(a,R)$, which
are pointwise limits of Riemann sums in the coordinates of $w$.

Since $\widetilde X^{(\al)}(t)=\Phi(\theta_t\ep)$ and
$\theta_h\ep\stackrel{d}{=}\ep$ for every $h\in\RR$, the process
$\widetilde X^{(\al)}$ is strictly stationary.
Fubini's theorem, $\EE[c_{\ep(s)}]=\kappa$ \textup(Lemma \ref{lem:cov}\textup)
and \eqref{eq:Kzeromass} give \eqref{eq:Xstat-mean}.
Finally, by \eqref{eq:centered-c} ($c_{\ep(t)}=\bar c+c\eta(t)$)
 and zero-mass equality \eqref{eq:Kzeromass} once more, the
constant part of the velocity contributes $\bar c\int_0^\infty K_\al(u)\rmd u=0$, so that
\begin{equation}\label{eq:Xstat-centered}
\widetilde X^{(\al)}(t)=c\int_0^\infty K_\al(u)\,\eta(t-u)\rmd u,
\qquad\text{whence}\qquad
\bigl|\widetilde X^{(\al)}(t)\bigr|\le c\,\|K_\al\|_{L^1(0, \infty)}=2cK_{1+\al}(t^*)
\end{equation}
by \eqref{eq:KL1}.

\textup{(ii)}
By the representation \eqref{eq:repr} and \eqref{eq:Kzeromass},
\[
X^{(\al)}(t)-\widetilde X^{(\al)}(t)
=x\rme^{-\beta t}-\int_t^\infty K_\al(u)\,c_{\ep(t-u)}\rmd u ,
\]
and for $t\ge t^*,$ by \eqref{eq:Katail}, the last term is at most
$M\int_t^\infty|K_\al(u)|\rmd u=MK_{1+\al}(t)$ in absolute value, which decays
at the rate given by \eqref{eq:Kasym+1}.

\textup{(iii)}
Let $t\ge2t^*$ and $h\ge0$. Splitting the integral in \eqref{eq:repr},
\[\aligned
X^{(\al)}(t+h)=&x\rme^{-\beta(t+h)}+Y_t(h)+R_t(h),\\
Y_t(h):=&\int_0^{t/2+h}K_\al(u)\,c_{\ep(t+h-u)}\rmd u ,\quad
R_t(h):=\int_{t/2+h}^{t+h}K_\al(u)\,c_{\ep(t+h-u)}\rmd u.
\endaligned\]
By \eqref{eq:Katail} and Theorem \ref{theo:kernel}\,\textup{(v)}, we have
 \[|R_t(h)|\le M\int_{t/2+h}^\infty|K_\al(u)|\rmd u
 =MK_{1+\al}(t/2+h)
 \le MK_{1+\al}(t/2).\]

In the same way,
\[\aligned
\widetilde X^{(\al)}(t+h)= &\widetilde Y_t(h)+\widetilde R_t(h),\\
\widetilde R_t(h)=&\int_{t/2+h}^\infty K_\al(u)c_{\ep(t+h-u)}\rmd u,\qquad |\widetilde R_t(h)|\le MK_{1+\al}(t/2).
\endaligned\]
Here $\widetilde Y_t(h)$ is defined in the same way from the stationary version.
$Y_t(h)$ and $\widetilde Y_t(h)$ are the same measurable functional of the
path $\{\ep(t/2+s)\}_{s\ge0}$, whose law is $\sum_j\pi_{ij}(t/2)\,\PP_j$
under $\PP_i$, by the Markov property, and $\sum_jp_j^*\,\PP_j$ under the
stationary version. By \eqref{eq:pi}, these two laws are at total variation distance
$p_{1-i}^*\rme^{-\la t}$.

Write $\mathbf X_t$ and $\widetilde{\mathbf X}_t$ for the two vectors with
components $X^{(\al)}(t+h_k)$ and $\widetilde X^{(\al)}(t+h_k)$,
$k=1,\dots,n$.
For any bounded function $\varphi\colon\RR^n\to\RR$ that  is Lipschitz continuous with
constant $L$ with respect to the $\ell^1$-norm, that is, satisfying 
\[
|\varphi(\bm x)-\varphi(\bm x')|\le L\sum_{k=1}^n|x_k-x_k'|,\qquad \bm x, \bm x'\in\RR^n,
\]
the three estimates yield
\[
\bigl|\EE_i\,\varphi(\mathbf X_t)-\EE\,\varphi(\widetilde{\mathbf X}_t)\bigr|
\le nL\cdot\bigl(|x|\rme^{-\beta t}+2MK_{1+\al}(t/2)\bigr)
+2\|\varphi\|_\infty\,p_{1-i}^*\rme^{-\la t}\ \longrightarrow\ 0,\quad \text{as $t\to\infty$}.
\]
Since $\widetilde{\mathbf X}_t\stackrel{d}{=}\widetilde{\mathbf X}_0$ by
\textup{(i)}, this is \eqref{eq:limlaw}.

Both processes are bounded uniformly in $t$, by Theorem \ref{theo:band} and
\eqref{eq:Xstat-band}, so the mixed moments converge as well.
\end{proof}

The constant in \eqref{eq:Xstat-band} cannot be improved\textup: the essential
supremum and infimum of $\widetilde X^{(\al)}(t)$ are $\pm2cK_{1+\al}(t^*)$.
Indeed, take $x=0$ in \textup{(ii)}. By Proposition \ref{prop:support}, for
every $\delta>0$ the event \[\{X^{(\al)}(t)>M_\al(t)-\delta\}\] has positive
probability, and on it
$\widetilde X^{(\al)}(t)>M_\al(t)-\delta-MK_{1+\al}(t)$ by \eqref{eq:coupling}.
The law of $\widetilde X^{(\al)}(t)$ does not depend on $t$, and
$M_\al(t)\to2cK_{1+\al}(t^*)$ by Theorem \ref{theo:band}\,\textup{(iii)}.
Letting $t\to\infty$ and then $\delta\downarrow0$ gives the supremum; the
infimum is treated symmetrically.

\begin{remark}\label{rem:stat-alpha1}
For $\al=1,$ the same construction yields a stationary
Kac--Ornstein--Uhlenbeck process
$\widetilde X^{(1)}(t)=\int_{-\infty}^t\rme^{-\beta(t-s)}c_{\ep(s)}\rmd s$\textup;
its mean is $\kappa/\beta$ (rather than $0$), it takes values within the attractor
$[\rho_1,\rho_0]$ described in Corollary \ref{cor:attractor}, and the convergence in
\textup{(ii)} is exponential (rather than $O(t^{\al-1})$).
This difference traces back to the factor $\tau^{1-\al}$ in \eqref{eq:LTKal},
which vanishes at $\tau=0$ only if $\al<1$.
\end{remark}

\begin{proposition}
\label{prop:acov-stat}
Let $0<\al<1$ and $\beta>0$. The autocovariance
$r(h):=\Cov\bigl(\widetilde X^{(\al)}(0),\widetilde X^{(\al)}(h)\bigr)$\textup,
$h\in\RR$\textup, of the stationary version is an even continuous
function\textup, and
\begin{equation}\label{eq:acov-double}
r(h)=\sigma^2\int_0^\infty\!\!\int_0^\infty K_\al(u)K_\al(v)\,
\rme^{-2\la|h+u-v|}\rmd u\rmd v .
\end{equation}
For $h\ge0$\textup,
\begin{equation}\label{eq:acov-stat}
r(h)=\sigma^2\int_0^\infty
\Bigl[K_\al(u)\,J_1(h+u)+K_\al(h+u)\,J_1(u)\Bigr]\rmd u\textup;
\end{equation}
in particular\textup, $r(0)=2\sigma^2\int_0^\infty K_\al(u)J_1(u)\rmd u$\textup, which is the
limit of the variances in Corollary \textup{\ref{cor:moments-limit}}.
\end{proposition}

\begin{proof}
By \eqref{def:Xstat} and Fubini's theorem,
\[
r(h)=\int_0^\infty\!\!\int_0^\infty K_\al(u)K_\al(v)\,
\Cov\bigl(c_{\ep(-u)},c_{\ep(h-v)}\bigr)\rmd u\rmd v,
\]
which, by the stationary identity \eqref{eq:cov-stat}, gives \eqref{eq:acov-double}.
Interchanging $u$ and $v$ shows that $r$ is even, and $r$ is continuous by
dominated convergence, the integrand being dominated by $|K_\al(u)K_\al(v)|$.

Let $h\ge0$. By \eqref{eq:acov-double},
\[
r(h)=\sigma^2\int_0^\infty K_\al(u)\left(\int_0^\infty K_\al(v)\,
\rme^{-2\la|h+u-v|}\rmd v \right)\rmd u,
\]
and we split the inner integral at $v=h+u$.
By \eqref{eq:I1ode} the part over $v<h+u$ is $J_1(h+u)$, and by
\eqref{def:Theta} the part over $v>h+u$ is $\Theta(h+u)$, so that
$r(h)=\sigma^2\int_0^\infty K_\al(u)\bigl[J_1(h+u)+\Theta(h+u)\bigr]\rmd u$.
Expanding $\Theta$ by \eqref{def:Theta} and substituting $v=u+w$ in the
resulting double integral,
\begin{align*}
\int_0^\infty\!\!\int_0^\infty\rme^{-2\la w}K_\al(u)K_\al(h+u+w)\rmd w\rmd u
&=\int_0^\infty\!\!\int_0^{v}\rme^{-2\la w}K_\al(v-w)K_\al(h+v)\rmd w\rmd v\\
&=\int_0^\infty K_\al(h+v)J_1(v)\rmd v ,
\end{align*}
which gives \eqref{eq:acov-stat}.
\end{proof}

By Proposition \ref{prop:stationary}\,\textup{(iii)}, $r$ is also the limit of
the autocovariance of Theorem \ref{theo:acov}\textup: for every $h\ge0$,
$x\in\RR$ and $i\in\{0,1\}$,
\[
\Cov_i\bigl(X^{(\al)}(t),X^{(\al)}(t+h)\bigr)\longrightarrow r(h),
\qquad t\to\infty .
\]

\subsection{Antipersistence}\label{sec:antipersistence}

A stationary process with integrable autocovariance $r$ is called
\emph{antipersistent} if $\int_{-\infty}^\infty r(h)\rmd h=0$, that is, if its
spectral density vanishes at the origin\textup; see \citet[Sect.~1.3]{LM}. The
positive correlations at short lags are then exactly cancelled by negative
correlations at long lags.

The spectral density $f_\al$ of $\widetilde X^{(\al)}$ is defined by
$r(h)=\int_{-\infty}^\infty \rme^{\rmi\omega h}f_\al(\omega)\rmd \omega.$

\begin{theorem}[Antipersistence]\label{theo:antipersistence}
Let $0<\al<1$ and $\beta>0$\textup, and let $\widetilde X^{(\al)}$ be defined by
\eqref{def:Xstat}.
\begin{enumerate}
\item[\textup{(i)}] $\widetilde X^{(\al)}$ has the spectral density
\begin{equation}\label{eq:specdens}
f_\al(\omega)=\frac{\sigma^2}{2\pi}\cdot\frac{4\la}{4\la^2+\omega^2}\cdot
\frac{|\omega|^{2-2\al}}{\omega^2+\beta^2},
\qquad \omega\in\RR;
\end{equation}
in particular
$f_\al(\omega)\sim\sigma^2|\omega|^{2-2\al}/(2\pi\la\beta^2)$ as $\omega\to0;$
\item[\textup{(ii)}] $r(0)>0;$
\item[\textup{(iii)}] $r\in L^1(\RR)$ and
\begin{equation}\label{eq:zerointegral}
\int_{-\infty}^{\infty}r(h)\rmd h=0 .
\end{equation}
\end{enumerate}
\end{theorem}

\begin{proof}
\textup{(i)}
Since $\sigma^2\rme^{-2\la|x|}
=\int_{-\infty}^\infty\rme^{\rmi\omega x}g(\omega)\rmd\omega$ with
$g(\omega):=\frac{\sigma^2}{2\pi}\frac{4\la}{4\la^2+\omega^2}$, substituting this
into \eqref{eq:acov-double} gives
\[
r(h)=\int_0^\infty\int_0^\infty K_\al(u)K_\al(v)
\left(\int_{-\infty}^\infty\rme^{\rmi\omega|h+u-v|}g(\omega)\rmd\omega\right)
\rmd u\rmd v,
\]
and, by Fubini's theorem \textup($K_\al$ and $g$ being integrable\textup),
\[
r(h)=\int_{-\infty}^\infty\rme^{\rmi\omega h}g(\omega)\,
\bigl|\widehat{K_\al}(\rmi\omega)\bigr|^2\rmd\omega,
\qquad
\widehat{K_\al}(\rmi\omega):=\int_0^\infty\rme^{-\rmi\omega v}K_\al(v)\rmd v .
\]
As $K_\al\in L^1(0,\infty)$, dominated convergence shows that
$\widehat{K_\al}(\rmi\omega)$ is the limit of $\widehat{K_\al}(\tau)$ as
$\tau\to\rmi\omega$ with $\RRe\,\tau>0$, so that by \eqref{eq:LTKal}
$\widehat{K_\al}(\rmi\omega)=(\rmi\omega)^{1-\al}/(\rmi\omega+\beta)$ and
$|\widehat{K_\al}(\rmi\omega)|^2=|\omega|^{2-2\al}/(\omega^2+\beta^2)$.
This is \eqref{eq:specdens}.

\textup{(ii)}
By \textup{(i)}, $r(0)=\int_{-\infty}^\infty f_\al(\omega)\rmd\omega>0$. 

\textup{(iii)}
By Tonelli's theorem and \eqref{eq:acov-double},
$\int_{-\infty}^\infty|r(h)|\rmd h\le\|K_\al\|_{L^1}^2\cdot\sigma^2/\la<\infty$\textup;
and by Fubini's theorem, integrating first in $h$,
\[
\int_{-\infty}^\infty r(h)\rmd h
=\frac{\sigma^2}{\la}\Bigl(\int_0^\infty K_\al(u)\rmd u\Bigr)^2=0;
\]
see \eqref{eq:Kzeromass}.
\end{proof}
\begin{remark}
Since $r$ is continuous and $r(0)>0$ by \textup{(ii)}, $r$ is positive in a
neighbourhood of the origin; as $\int_{-\infty}^\infty r(h)\rmd h=0$ by
\textup{(iii)}, $r$ must take negative values somewhere.
\end{remark}

\begin{remark}[Bounded variation and antipersistence]\label{rem:antipersistence}
On compact intervals,
the trajectories of the fKOU process have bounded variation: by
\eqref{eq:path}, each trajectory is a finite linear combination of
$\rme^{-\beta t}$ and of shifts $K_{1+\al}(\cdot-\tau_k)$, each of which is
unimodal. Nevertheless, its stationary correlation structure is that of an
antipersistent model.

In the case $\tfrac12<\al<1,$ the process is antipersistent also in the stricter sense
of \citet[Definition~1.4]{LM}, with a memory parameter
$d=\al-1\in(-\tfrac12,0)$. This means that its low-frequency behaviour is that of
a fractional Gaussian noise with Hurst exponent $H=d+\tfrac12=\al-\tfrac12$
\citep{MvN}. 

For $0<\al\le\tfrac12$ the spectral density vanishes at the origin faster
than this definition allows.
When $\al=1$, by contrast, $r$ decays exponentially and
$\int_{-\infty}^\infty r(h)\rmd h=\sigma^2/(\la\beta^2)>0$; the behaviour of
the process is therefore discontinuous at $\al=1$, a consequence of
\eqref{eq:Kzeromass}.
\end{remark}

The stationary variance can be computed in closed form from the spectral
density \eqref{eq:specdens}.

\begin{proposition}[Stationary variance]\label{prop:r0}
Let $0<\al<1$ and $\beta>0$. Then
\begin{equation}\label{eq:r0closed}
r(0)=\Var\widetilde X^{(\al)}(0)
=\frac{2\la\,\sigma^2}{4\la^2-\beta^2}\cdot
\frac{(2\la)^{1-2\al}-\beta^{1-2\al}}{\cos(\pi\al)} ,
\end{equation}
the singularities at $\al=\tfrac12$ and at $2\la=\beta$ are removable\textup:
there the right-hand side is to be read as
\[
\frac{2\la\,\sigma^2}{4\la^2-\beta^2}\cdot\frac2\pi\log\frac{2\la}{\beta}
\qquad\text{and}\qquad
\frac{(1-2\al)\,\sigma^2}{2\beta^{2\al}\cos(\pi\al)} ,
\]
respectively\textup, and as $\sigma^2/(\pi\beta)$ if both conditions hold.
\end{proposition}
\begin{proof}
By Theorem \ref{theo:antipersistence}\,\textup{(i)},
$r(0)=\int_{-\infty}^\infty f_\al(\omega)\rmd\omega$, which is a continuous function of
$(\al,\beta,\la)$ by dominated convergence. Consequently, it suffices to prove
\eqref{eq:r0closed} for $\al\ne\tfrac12$ and $2\la\ne\beta$\textup; the values
at the exceptional parameters are then the limits of the right-hand side, and
l'H\^opital's rule gives the three expressions in the statement.

Put $p:=3-2\al$. By partial fractions,
\[
r(0)
=\frac{4\la\,\sigma^2}{\pi\,(4\la^2-\beta^2)}\int_0^\infty\omega^{p-1}
\Bigl[\frac{1}{\omega^2+\beta^2}-\frac{1}{\omega^2+4\la^2}\Bigr]\rmd\omega .
\]
For $0<p<2$ the two terms may be integrated separately, and
$\int_0^\infty\frac{\omega^{p-1}}{\omega^2+a^2}\rmd\omega
=\frac{\pi a^{p-2}}{2\sin(\pi p/2)}$ gives the value
$\frac{\pi}{2\sin(\pi p/2)}\bigl(\beta^{p-2}-(2\la)^{p-2}\bigr)$ for the
integral. The integral of the difference converges for every $0<p<4$, the
integrand being $O(\omega^{p-5})$ at infinity, and is analytic in $p$ there\textup;
by analytic continuation the same value holds for $0<p<4$, $p\ne2$, a range
that covers every $\al\in(0,1]$.
Since $p-2=1-2\al$ and $\sin(\pi p/2)=-\cos(\pi\al)$, this is
\eqref{eq:r0closed}.
\end{proof}

\begin{remark}
For $\al=1$ formula \eqref{eq:r0closed} gives
$\sigma^2/\bigl(\beta(\beta+2\la)\bigr)$\textup, the stationary variance of the
Kac--Ornstein--Uhlenbeck process.
\end{remark}

\subsection{A law of large numbers}\label{sec:lln}

By Proposition \ref{prop:stationary} and Remark \ref{rem:stat-alpha1}, for $\beta>0$ the process $X^{(\al)}$ does
not grow\textup: it settles into the stationary regime.
The driving process
$\sF^{(\al)}$ does grow, and obeys a law of large numbers in which the
parameters are not rescaled.

\begin{proposition}[Law of large numbers]\label{prop:lln}
Let $0<\al\le1$. For each initial state\textup, almost surely\textup,
\begin{equation}\label{eq:lln}
\lim_{t\to\infty}\frac{\sF^{(\al)}(t)}{t^\al}=\frac{\kappa}{\Gamma(1+\al)} .
\end{equation}
Consequently\textup, for $\beta>0$ and every $x\in\RR$\textup, almost
surely\textup,
\begin{equation}\label{eq:lln-X}
\lim_{t\to\infty}\frac1{t^\al}\int_0^tX^{(\al)}(s)\rmd s
=\frac{\kappa}{\beta\,\Gamma(1+\al)} .
\end{equation}
\end{proposition}

\begin{proof}
Write $A(t):=\int_0^tc_{\ep(s)}\rmd s$, so that $|A(t)-A(s)|\le M|t-s|$ and
$A(t)/t\to\kappa$ almost surely (see \eqref{eq:ergodic}).
Integrating by parts in the definition \eqref{eq:fKOU} of $\sF^{(\al)}$ and
substituting $s=tv$,
\[
\frac{\Gamma(\al)\,\sF^{(\al)}(t)}{t^\al}
=\frac{A(t)}{t}+(1-\al)\int_0^1(1-v)^{\al-2}\,\frac{A(t)-A(tv)}{t}\,\rmd v .
\]

On the almost sure event where $A(u)/u\to\kappa$ as $u\to\infty$,
the integrand converges to
$\kappa(1-v)^{\al-1}$ for every $v\in(0,1)$ and is dominated by
$M(1-v)^{\al-1}$. By dominated convergence the right-hand side tends to
$\kappa+(1-\al)\kappa/\al=\kappa/\al$, which is \eqref{eq:lln} since
$\al\Gamma(\al)=\Gamma(1+\al)$.

For the second statement, \eqref{eq:fKOU} reads
$\beta\int_0^tX^{(\al)}(s)\rmd s=x+\sF^{(\al)}(t)-X^{(\al)}(t)$, and
$X^{(\al)}$ is bounded by Theorem \ref{theo:band}\textup; dividing by $t^\al$
and using \eqref{eq:lln} gives \eqref{eq:lln-X}.
\end{proof}

In the symmetric case, $c_0=-c_1=c$ and
$\la_0=\la_1=\la$,
one has $\kappa=0$, and
\eqref{eq:lln} says only that $\sF^{(\al)}(t)=o(t^\al)$. The fluctuations
are the subject of the next section.

\section{Kac's scaling limit}\label{sec:kac}
\setcounter{equation}{0}

Throughout this section the velocity process is symmetric, $c_0=-c_1=c$ and
$\la_0=\la_1=\la$, and we study the fKOU process under \emph{Kac's scaling}
\begin{equation}\label{def:scaling}
c,\la\to\infty\quad\text{with}\quad c^2/\la\to1 .
\end{equation}
The initial law of $\ep$ is arbitrary and fixed, and in this section $\PP$ and
$\EE$ refer to it. We also fix $T>0$, $\beta\ge0$ and $x\in\RR$.

In the symmetric case, \eqref{eq:centered-c} reads 
$c_{\ep(t)}=c\,\eta(t),$ where $\eta(t)=1-2\ep(t)\in\{-1, +1\}$, so
that the integrated telegraph process \eqref{def:T} is
$\sT(t)=c\int_0^t\eta(u)\rmd u$. Under \eqref{def:scaling} it converges to a
standard Brownian motion $W$. This goes back to \citet{Kac}\textup; see also
\citet{MPRF,RK}.

Our goal is to transfer this convergence to $X^{(\al)}$. By the
representation \eqref{eq:repr}, $X^{(\al)}(t)-x\rme^{-\beta t}$ is the integral
of $K_\al(t-\cdot)$ against $\sT$, and an integration by parts turns it into
$\hat h[\sT](t)$, where $\hat h$ is a linear operator that involves only the values of
its argument and is bounded on $\gamma$-H\"older paths as soon as
$\gamma>1-\al$ (Section \ref{sec:kac-operator}). The paths of $W$ are
$\gamma$-H\"older only for $\gamma<\tfrac12$. We therefore first prove that $\sT$
converges to $W$ in H\"older norm of every order $\gamma<\tfrac12$ (Section
\ref{sec:kac-telegraph})\textup; the continuous mapping theorem then gives the
limit of $X^{(\al)}$ whenever the window $1-\al<\gamma<\tfrac12$ is non-empty,
that is, for $\al>\tfrac12$ (Section \ref{sec:kac-limit}). For $\al\le\tfrac12$
the window is empty, and Section \ref{sec:kac-small} shows that this is not a
defect of the method.

For $0<\gamma\le1$ we write
\[
[f]_\gamma:=\sup_{0\le s<t\le T}\frac{|f(t)-f(s)|}{(t-s)^\gamma}
\]
for the H\"older seminorm of $f\colon[0,T]\to\RR$, and $C^\gamma[0,T]$ for the
space of functions with $[f]_\gamma<\infty$, a Banach space for the norm
$|f(0)|+[f]_\gamma$. For $\gamma<1$ we denote by $C^{0,\gamma}[0,T]$ its closed
subspace of functions such that
\[
\lim_{\delta\downarrow0}\ \sup_{0<t-s\le\delta}\frac{|f(t)-f(s)|}{(t-s)^\gamma}=0 .
\]
Unlike $C^\gamma[0,T]$, the space $C^{0,\gamma}[0,T]$ is separable, which is why
weak convergence is stated there\textup; it contains $C^{\gamma'}[0,T]$ for
every $\gamma'>\gamma$.

\subsection{The telegraph process in H\"older norm}\label{sec:kac-telegraph}

We prove that $\sT$ converges to $W$ in $C^{0,\gamma}[0,T]$ for every
$\gamma<\tfrac12$. The proof has the two usual parts: a moment bound, which
gives tightness (Lemmas \ref{lem:Tmoments} and \ref{lem:tight}), and the
convergence of the finite-dimensional distributions, which we read off the
Cattaneo equation satisfied by the characteristic functions.

\begin{lemma}[Moments of the integrated telegraph process]\label{lem:Tmoments}
Let $m\ge1$ be an integer.
\begin{enumerate}
\item[\textup{(i)}] For all $0\le u_1\le\dots\le u_{2m}$\textup,
\begin{equation}\label{eq:sigcorr}
\EE\Bigl[\prod_{i=1}^{2m}\eta(u_i)\Bigr]
=\prod_{j=1}^{m}\rme^{-2\la(u_{2j}-u_{2j-1})},
\end{equation}
whatever the initial law of $\ep$.
\item[\textup{(ii)}] For all $0\le s\le t$\textup,
\begin{equation}\label{eq:Tmoments}
\EE\bigl|\sT(t)-\sT(s)\bigr|^{2m}
\ \le\ (2m-1)!!\ \Bigl(\frac{c^2}{\la}\Bigr)^{m}(t-s)^{m} .
\end{equation}
\end{enumerate}
\end{lemma}

\begin{proof}
\textup{(i)}
Since $\la_0=\la_1=\la$, the switching times of $\ep$ form a Poisson process of
rate $\la$ that is independent of $\ep(0)$. Let $N(a,b]$ be the number of
switchings in $(a,b]$\textup; then $\eta(a)\eta(b)=(-1)^{N(a,b]}$ and
$\EE\bigl[(-1)^{N(a,b]}\bigr]=\rme^{-2\la(b-a)}$. Grouping the $2m$ factors in
consecutive pairs,
\[
\prod_{i=1}^{2m}\eta(u_i)=\prod_{j=1}^m\eta(u_{2j-1})\eta(u_{2j})
=\prod_{j=1}^m(-1)^{N(u_{2j-1},u_{2j}]} ,
\]
a product over disjoint intervals, hence of independent factors. This gives
\eqref{eq:sigcorr}, and the initial law plays no role.

\textup{(ii)}
Expanding the $2m$-th power and using the symmetry of the integrand under
permutations,
\begin{align*}
\EE\bigl(\sT(t)-\sT(s)\bigr)^{2m}
&=c^{2m}\int_{[s,t]^{2m}}\EE\Bigl[\prod_{i=1}^{2m}\eta(u_i)\Bigr]\rmd\bm u\\
&=(2m)!\,c^{2m}\int_{s\le u_1\le\dots\le u_{2m}\le t}\
\prod_{j=1}^m\rme^{-2\la(u_{2j}-u_{2j-1})}\rmd\bm u .
\end{align*}
Pass to the variables $u_1\le u_3\le\dots\le u_{2m-1}$ and the $m$ gaps
$a_j:=u_{2j}-u_{2j-1}\ge0$. Dropping the constraints $u_{2j}\le u_{2j+1}$ and
$u_{2m}\le t$ enlarges the domain of integration, and the integrand is
positive, so the last integral is at most the product of the volume
$(t-s)^m/m!$ of $\{s\le u_1\le u_3\le\dots\le u_{2m-1}\le t\}$ and of
$\bigl(\int_0^\infty\rme^{-2\la a}\rmd a\bigr)^m=(2\la)^{-m}$. Hence
\[
\EE\bigl(\sT(t)-\sT(s)\bigr)^{2m}
\le\frac{(2m)!}{m!\,2^m}\Bigl(\frac{c^2}{\la}\Bigr)^m(t-s)^m
=(2m-1)!!\,\Bigl(\frac{c^2}{\la}\Bigr)^m(t-s)^m. \qedhere
\]
\end{proof}

For $c^2/\la=1$ the right-hand side of \eqref{eq:Tmoments} is exactly
$\EE|W(t)-W(s)|^{2m}$\textup: the whole of Kac's scaling is visible in this one
bound.

A bound of this kind gives tightness in H\"older norm. The criterion goes back
to \citet{Lamperti}\textup; see also \citet[Corollary~1]{Hamadouche}. We include the
short proof because it yields the explicit estimate \eqref{eq:holdermom}.

\begin{lemma}[Tightness in H\"older norm]\label{lem:tight}
Let $\sX$ be a family of processes with continuous paths on $[0,T]$ such that
$X(0)=0$ for every $X\in\sX$. Suppose that\textup, for some integer $m\ge2$ and
some $C<\infty$\textup,
\begin{equation}\label{eq:kolm}
\EE\bigl|X(t)-X(s)\bigr|^{2m}\le C|t-s|^m,
\qquad 0\le s,t\le T,\quad X\in\sX .
\end{equation}
Then\textup, for every $0<\gamma<\frac{m-1}{2m}$\textup:
\begin{enumerate}
\item[\textup{(i)}] $\sup_{X\in\sX}\EE\,[X]_\gamma^{2m}<\infty$\textup;
\item[\textup{(ii)}] the laws of the processes $X\in\sX$ are tight in
$C^{0,\gamma}[0,T]$.
\end{enumerate}
\end{lemma}

\begin{proof}
\textup{(i)}
Let $\mathbb D_n:=\{kT2^{-n}:0\le k\le2^n\}$ be a partition of $[0, T]$, and let
$\mathbb D:=\bigcup_n\mathbb D_n$. For $X\in\sX$ put
\[
\Delta_n:=\max_{1\le k\le2^n}\bigl|X(kT2^{-n})-X((k-1)T2^{-n})\bigr| .
\]
By \eqref{eq:kolm},
\[
\EE\bigl[\Delta_n^{2m}\bigr]
\le\sum_{k=1}^{2^n}\EE\bigl|X(kT2^{-n})-X((k-1)T2^{-n})\bigr|^{2m}
\le2^n\,C(T2^{-n})^m=CT^m2^{-n(m-1)},
\]
that is,
\begin{equation}\label{eq:dyadic}
\|\Delta_n\|_{L^{2m}(\PP)}\le\bigl(CT^m\bigr)^{1/(2m)}\,2^{-n\frac{m-1}{2m}} .
\end{equation}
Let $s<t$ be points of $\mathbb D$ with $t-s\le T2^{-N}$, let $s'$ be the
smallest point of $\mathbb D_N$ with $s'\ge s$, and let $t'$ be the largest
point of $\mathbb D_N$ with $t'\le t$. Halving repeatedly, $s$ is joined to
$s'$ by at most one step of each level $n>N$, so that
$|X(s')-X(s)|\le\sum_{n>N}\Delta_n$, and likewise for $t$ and $t'$\textup;
moreover $s'$ and $t'$ coincide or are adjacent in $\mathbb D_N$, so that
$|X(s')-X(t')|\le\Delta_N$. Hence
\begin{equation}\label{eq:chaining}
|X(t)-X(s)|\ \le\ 2\sum_{n\ge N}\Delta_n .
\end{equation}
Given $s\ne t$ in $\mathbb D$, choose $N\ge0$ with
$T2^{-N-1}<|t-s|\le T2^{-N}$. By \eqref{eq:chaining}, and since
$2^{N\gamma}\le2^{n\gamma}$ for $n\ge N$,
\[
\frac{|X(t)-X(s)|}{|t-s|^\gamma}
\le\frac{2\sum_{n\ge N}\Delta_n}{(T2^{-N-1})^\gamma}
=2^{1+\gamma}T^{-\gamma}\,2^{N\gamma}\!\sum_{n\ge N}\Delta_n
\ \le\ 2^{1+\gamma}T^{-\gamma}\sum_{n\ge0}2^{n\gamma}\Delta_n .
\]
The paths being continuous, the supremum of the left-hand side over
$\mathbb D$ is $[X]_\gamma$, so by Minkowski's inequality and
\eqref{eq:dyadic}
\begin{equation}\label{eq:holdermom}
\bigl\|[X]_\gamma\bigr\|_{L^{2m}(\PP)}
\le2^{1+\gamma}T^{-\gamma}\bigl(CT^m\bigr)^{1/(2m)}
\sum_{n\ge0}2^{n\left(\gamma-\frac{m-1}{2m}\right)}=:R .
\end{equation}
The series converges because $\gamma<\frac{m-1}{2m}$, and $R$ does not depend
on $X$.

\textup{(ii)}
Fix $\gamma<\gamma'<\frac{m-1}{2m}$ and put
$A_M:=\{f:f(0)=0,\ [f]_{\gamma'}\le M\}$. If $f_k\in A_M$, then $\{f_k\}$ is
bounded and equicontinuous, so by the Arzel\`a--Ascoli theorem a subsequence
converges uniformly to some $f$, necessarily in $A_M$\textup; and then
\begin{align*}
\frac{|(f_k-f)(t)-(f_k-f)(s)|}{|t-s|^\gamma}
&\le\Bigl(\frac{|(f_k-f)(t)-(f_k-f)(s)|}{|t-s|^{\gamma'}}\Bigr)^{\gamma/\gamma'}
\bigl|(f_k-f)(t)-(f_k-f)(s)\bigr|^{1-\gamma/\gamma'}\\
&\le(2M)^{\gamma/\gamma'}\bigl(2\|f_k-f\|_\infty\bigr)^{1-\gamma/\gamma'}
\ \longrightarrow\ 0 ,
\end{align*}
so the convergence holds in $C^\gamma[0,T]$. 

Moreover
$\sup_{0<t-s\le\delta}|f(t)-f(s)|/(t-s)^\gamma\le M\delta^{\gamma'-\gamma}\to0$
for $f\in A_M$, so $A_M\subset C^{0,\gamma}[0,T]$, and $A_M$ is compact there.

By \eqref{eq:holdermom} applied with $\gamma'$ and Chebyshev's inequality,
$\sup_{X\in\sX}\PP\{X\notin A_M\}\le R^{2m}M^{-2m}$, where $R$ is the constant
of \eqref{eq:holdermom} for $\gamma'$, and tightness follows.
\end{proof}

We can now state Kac's theorem in the form we need. The convergence of
one-dimensional distributions is \citet[Theorem~2.7]{RK} and its functional
version in $C[0,T]$ is \citet[Remark~2.11]{RK}; we need it in the stronger
H\"older norm.

\begin{theorem}[Kac's limit in H\"older norm]\label{theo:kac}
Under the scaling \eqref{def:scaling}\textup, for every
$0<\gamma<\tfrac12$\textup,
\[
\sT\Longrightarrow W\qquad\text{in }C^{0,\gamma}[0,T],
\]
where $W$ is a standard Brownian motion\textup; in particular
$\sT\Longrightarrow W$ in $C[0,T]$.
\end{theorem}

\begin{proof}
\emph{Step 1\textup: one-time characteristic functions.}
For $k,j\in\{0,1\}$ and $\xi\in\RR$ put
\[
\mu_{kj}(t,\xi):=\EE_j\bigl[\rme^{\rmi\xi\sT(t)}\1_{\{\ep(t)=k\}}\bigr],
\qquad \mathbf M(t,\xi):=\bigl(\mu_{kj}(t,\xi)\bigr)_{k,j} ,
\]
where $\EE_j$ is the expectation given $\ep(0)=j$.
Conditioning on the first switching from the state $j$,
\[
\mu_{kj}(t,\xi)=\rme^{-\la t}\rme^{\rmi\xi c_jt}\delta_k(j)
+\int_0^t\la\rme^{-\la s}\rme^{\rmi\xi c_js}\,\mu_{k,1-j}(t-s,\xi)\rmd s ,
\]
whence, differentiating in $t$,
\begin{equation}\label{eq:kacsys}
\partial_t\mu_{kj}=\rmi\xi c_j\mu_{kj}-\la\mu_{kj}+\la\mu_{k,1-j},
\qquad\mu_{kj}(0,\xi)=\delta_k(j),
\end{equation}
where $\delta_k(j)$ is the Kronecker delta.

Fix $k$ and set $S:=\mu_{k0}+\mu_{k1}$ and $D:=\mu_{k0}-\mu_{k1}$. Since
$c_0=-c_1=c$, \eqref{eq:kacsys} becomes
\begin{equation}\label{eq:SD}
S'=\rmi\xi cD,\qquad D'=\rmi\xi cS-2\la D,\qquad S(0)=1,\quad D(0)=(-1)^k ,
\end{equation}
and substituting the first equation into the second shows that $S$ solves the
\emph{Cattaneo equation}
\begin{equation}\label{eq:cattaneo}
S''+2\la S'+\xi^2c^2S=0,\qquad S(0)=1,\quad S'(0)=\rmi\xi c\,(-1)^k .
\end{equation}
Its characteristic roots are $r_\pm=-\la\pm\sqrt{\la^2-\xi^2c^2}$. Under
\eqref{def:scaling}, with $\xi$ fixed, $\xi^2c^2<\la^2$ for $\la$ large\textup;
the roots are then real and distinct,
\begin{equation}\label{eq:Ssol}
S(t)=\frac{\bigl(S'(0)-r_-\bigr)\rme^{r_+t}
+\bigl(r_+-S'(0)\bigr)\rme^{r_-t}}{r_+-r_-} ,
\end{equation}
and
\begin{equation}\label{eq:roots}
r_+=\frac{-\xi^2c^2}{\la+\sqrt{\la^2-\xi^2c^2}}
\longrightarrow-\frac{\xi^2}{2},
\qquad r_-\sim-2\la,\qquad r_+-r_-\sim2\la .
\end{equation}
Since $S'(0)=O(c)=O(\sqrt\la)$, the first coefficient in \eqref{eq:Ssol} tends
to $1$ and the second is $O(\la^{-1/2})$, while $|\rme^{r_-t}|\le1$\textup;
hence $S(t)\to\rme^{-t\xi^2/2}$ for every $t\ge0$.
For $D$, the second equation in \eqref{eq:SD} gives
\[
D(t)=(-1)^k\rme^{-2\la t}
+\rmi\xi c\int_0^t\rme^{-2\la(t-s)}S(s)\rmd s .
\]
As $|\mu_{kj}|\le1$ we have $|S|\le2$, so that
$|D(t)|\le\rme^{-2\la t}+|\xi|\,c/\la\to0$ for every $t>0$.
Since $\mu_{k0}=\tfrac12(S+D)$ and $\mu_{k1}=\tfrac12(S-D)$,
\begin{equation}\label{eq:mulim}
\mu_{kj}(t,\xi)\ \longrightarrow\ \tfrac12\rme^{-t\xi^2/2},
\qquad t>0,\quad\xi\in\RR,
\end{equation}
for all four pairs $(k,j)$.

\emph{Step 2\textup: finite-dimensional distributions.}
Let $0=t_0<t_1<\dots<t_d\le T$ and $\xi_1,\dots,\xi_d\in\RR$. The increment
$\sT(t_k)-\sT(t_{k-1})$ is a function of $\ep$ on $[t_{k-1},t_k]$ alone, so the
Markov property of $\ep$ gives the transfer-matrix identity
\begin{equation}\label{eq:transfer}
\EE\Bigl[\prod_{k=1}^d\rme^{\rmi\xi_k(\sT(t_k)-\sT(t_{k-1}))}\Bigr]
=\mathbf 1^\top \mathbf M(t_d-t_{d-1},\xi_d)\cdots
\mathbf M(t_1-t_0,\xi_1)\,\bm p,
\end{equation}
where $\mathbf 1:=(1,1)^\top$ and
$\bm p:=\bigl(\PP\{\ep(0)=0\},\PP\{\ep(0)=1\}\bigr)^\top$ is the initial law of
$\ep$.
By \eqref{eq:mulim}, $\mathbf M(t,\xi)\to\tfrac12\rme^{-t\xi^2/2}\mathbf J$,
where $\mathbf J$ is the $2\times2$ matrix of ones. Since
$\mathbf J\bm p=\mathbf 1$ for every probability vector $\bm p$ and
$\mathbf J\mathbf 1=2\cdot\mathbf 1$, the right-hand side of
\eqref{eq:transfer} tends to
$\prod_{k=1}^d\rme^{-(t_k-t_{k-1})\xi_k^2/2}$, the joint characteristic
function of the increments of $W$. So the finite-dimensional distributions of
$\sT$ converge to those of $W$.

\emph{Step 3\textup: tightness and identification of the limit.}
Fix $0<\gamma<\tfrac12$. The space $C^{0,\gamma}[0,T]$ being metrizable, it is
enough to argue along an arbitrary sequence $(c_n,\la_n)$ obeying
\eqref{def:scaling}, that is, $c_n,\la_n\to\infty$ and $c_n^2/\la_n\to1$.

 There is $n_0$ with $c_n^2/\la_n\le2$ for $n\ge n_0$, and
then \eqref{eq:Tmoments} gives \eqref{eq:kolm} with $C=(2m-1)!!\,2^m$, for
every $m$. Choose $m\ge2$ with $\gamma<\frac{m-1}{2m}$, which is possible
because $\frac{m-1}{2m}\uparrow\tfrac12$. The process $\sT$ has continuous
paths and $\sT(0)=0$, so by Lemma \ref{lem:tight}\,\textup{(ii)} the laws with
$n\ge n_0$ are tight in $C^{0,\gamma}[0,T]$\textup; each of the finitely many
others is tight by itself, that space being Polish.

By Prohorov's theorem \citep[Theorem~5.1]{Billingsley}, the sequence of laws of $\sT$
under $(c_n,\la_n)$ is relatively compact,
that is, every subsequence has a
further subsequence that converges weakly in $C^{0,\gamma}[0,T]$ to some law
$\mu$. Since the evaluation maps $f\mapsto f(t)$ are continuous on
$C^{0,\gamma}[0,T]$, the continuous mapping theorem and \emph{Step 2} show that
the finite-dimensional distributions of $\mu$ are those of $W$.
They determine $\mu$\textup: the space being separable, its Borel
$\sigma$-field is generated by the closed balls, and a closed ball
$\{f\colon|f(0)-g(0)|+[f-g]_\gamma\le r\}$ is described by the values of $f$ at
rational times.

On the other hand, the paths of $W$ belong to
$C^{\gamma'}[0,T]\subset C^{0,\gamma}[0,T]$ for $\gamma<\gamma'<\tfrac12$, so
the law of $W$ is a Borel probability measure on $C^{0,\gamma}[0,T]$ with the
same finite-dimensional distributions. Hence $\mu$ is the law of $W$, and
since the limit does not depend on the subsequence, the whole family
converges.
\end{proof}

\subsection{The memory kernel as an operator}\label{sec:kac-operator}

We now make precise the integration by parts announced at the beginning of the
section. The operator is defined on H\"older paths of exponent $\gamma>1-\al$, for
which the integral below converges absolutely; on Lipschitz paths it coincides
with the convolution of $K_\al$ with the derivative.

\begin{lemma}[The memory kernel as an operator on H\"older paths]\label{lem:hbdd}
Let $0<\al\le1$ and $1-\al<\gamma\le1$. For $g\in C^\gamma[0,T]$ put
$\hat h[g](0):=0$ and\textup, for $0<t\le T$\textup,
\begin{equation}\label{eq:byparts}
\hat h[g](t):=K_\al(t)\bigl(g(t)-g(0)\bigr)
+\int_0^tK_\al'(t-s)\bigl(g(s)-g(t)\bigr)\rmd s .
\end{equation}
\begin{enumerate}
\item[\textup{(i)}] The integral in \eqref{eq:byparts} converges absolutely\textup,
$\hat h[g]\in C[0,T]$\textup, and
$\|\hat h[g]\|_\infty\le\bigl(C_1+C_2\bigr)[g]_\gamma$\textup, where
\[
C_1:=\sup_{0<t\le T}|K_\al(t)|\,t^\gamma<\infty,
\qquad
C_2:=\int_0^T|K_\al'(u)|\,u^\gamma\rmd u<\infty .
\]
Thus $\hat h$ is a bounded linear operator from $C^\gamma[0,T]$ to $C[0,T]$.
\item[\textup{(ii)}] If $g$ is Lipschitz\textup, then
$\hat h[g](t)=\int_0^tK_\al(t-s)\,g'(s)\rmd s$ for $0\le t\le T$.
\end{enumerate}
\end{lemma}

\begin{proof}
\emph{Finiteness of $C_1$ and $C_2$.}
For $\al=1$ both $K_1(u)=\rme^{-\beta u}$ and $K_1'$ are bounded and everything
below is immediate, so let $0<\al<1$. By Theorem
\ref{theo:kernel}\,\textup{(ii)}, $K_\al(u)\sim u^{\al-1}/\Gamma(\al)$ as
$u\downarrow0$, and by \eqref{eq:Kder},
$K_\al'(u)=\frac{\al-1}{\Gamma(\al)}u^{\al-2}-\beta K_\al(u)=O(u^{\al-2})$.
Since $\al-1+\gamma>0$, it follows that $C_1$ and $C_2$ are finite, and also that 
\begin{equation}\label{eq:C1C2}
\sup_{0<u\le\zeta}|K_\al(u)|\,u^\gamma\longrightarrow0,
\qquad
\int_0^\zeta|K_\al'(u)|\,u^\gamma\rmd u\longrightarrow0
\qquad(\zeta\downarrow0).
\end{equation}

\textup{(i)}
Substituting $u=t-s$,
\[
\hat h[g](t)=K_\al(t)\bigl(g(t)-g(0)\bigr)
+\int_0^tK_\al'(u)\bigl(g(t-u)-g(t)\bigr)\rmd u .
\]
The absolute value of the first term does not exceed 
$[g]_\gamma|K_\al(t)|t^\gamma\le C_1[g]_\gamma.$ 
The integrand of the second term is bounded in absolute value by
$[g]_\gamma|K_\al'(u)|u^\gamma$ uniformly in $t$, which ensures absolute
convergence and an estimate of the form $C_2[g]_\gamma$. 

By \eqref{eq:C1C2}, both terms tend
to $0$ as $t\downarrow0$, so $\hat h[g]$ is continuous at the origin. If $t_n\to t>0$,
the first term is continuous at $t$, and
$\1_{\{u<t_n\}}K_\al'(u)\bigl(g(t_n-u)-g(t_n)\bigr)$ converges for every $u\ne t$
to the corresponding expression at $t$, under the same integrable bound\textup;
dominated convergence gives the continuity of the second term.

\textup{(ii)}
Let $g$ be Lipschitz and $0<t\le T$. For $0<\delta<t$, integrating by parts on
$[0,t-\delta]$,
\[
\int_0^{t-\delta}K_\al(t-s)\,g'(s)\rmd s
=\Bigl[K_\al(t-s)\bigl(g(s)-g(t)\bigr)\Bigr]_{s=0}^{s=t-\delta}
+\int_0^{t-\delta}K_\al'(t-s)\bigl(g(s)-g(t)\bigr)\rmd s .
\]
As $\delta\downarrow0$, the boundary term at $s=t-\delta$ is
$K_\al(\delta)\bigl(g(t-\delta)-g(t)\bigr)=O(\delta^{\al})\to0$ and the one at
$s=0$ gives $K_\al(t)\bigl(g(t)-g(0)\bigr)$\textup; the integral on the left
converges because $g'$ is bounded and $K_\al\in L^1(0,T)$, and the one on the
right by \textup{(i)}. This is \eqref{eq:byparts}.
\end{proof}

\begin{remark}[Sharpness of $\gamma>1-\al$]\label{rem:gamma-sharp}
For $0<\al<1$ the condition $\gamma>1-\al$ cannot be relaxed. At
$\gamma=1-\al$ one has $C_2=\infty$, and $\hat h[g]$ need not be continuous at the
origin\textup: for $\beta=0$ and $g(s)=s^{1-\al}$, which belongs to
$C^{1-\al}[0,T]$, the integral in \eqref{eq:byparts} converges for every $t>0$
and $\hat h[g](t)=\Gamma(2-\al)$ for all $t>0$, whereas $\hat h[g](0)=0$.
\end{remark}

\subsection{The scaling limit}\label{sec:kac-limit}

Theorem \ref{theo:kac} and Lemma \ref{lem:hbdd} combine into the main result of
this section.

\begin{theorem}[Scaling limit of the fKOU process]\label{theo:kac-limit}
Let $\tfrac12<\al\le1$\textup, $\beta\ge0$ and $x\in\RR$. In the symmetric
case\textup, under Kac's scaling \eqref{def:scaling} and for every initial law
of $\ep$\textup,
\begin{equation}\label{eq:kac-limit}
X^{(\al)}\Longrightarrow Y^{(\al)}\quad\text{in }C[0,T],
\qquad
Y^{(\al)}(t):=x\rme^{-\beta t}+\int_0^tK_\al(t-s)\rmd W(s),
\end{equation}
where $W$ is a standard Brownian motion and the stochastic convolution is taken
in its continuous version.
The process $Y^{(\al)}(t)-x\rme^{-\beta t},\;t\ge0,$ is centred Gaussian with
covariance
\[
\int_0^{t_1\wedge t_2}K_\al(t_1-s)\,K_\al(t_2-s)\rmd s .
\]
For $\al=1$\textup, $Y^{(1)}$ is the Ornstein--Uhlenbeck process. 

For
$\beta=0$\textup, $Y^{(\al)}-x$ is the Riemann--Liouville fractional Brownian
motion $\frac{1}{\Gamma(\al)}\int_0^t(t-s)^{\al-1}\rmd W(s)$ of
\citet{Levy}\textup, whose variance is
$t^{2\al-1}/\bigl((2\al-1)\Gamma(\al)^2\bigr)$.
\end{theorem}

\begin{proof}
\emph{Step 1\textup: reduction to the telegraph process.}
The process $\sT$ is Lipschitz, with $\sT'=c_{\ep(\cdot)}$ almost everywhere,
so \eqref{eq:repr} and Lemma \ref{lem:hbdd}\,\textup{(ii)} give
\begin{equation}\label{eq:XhT}
X^{(\al)}(t)-x\rme^{-\beta t}=\hat h\bigl[\sT\bigr](t),\qquad 0\le t\le T .
\end{equation}
Since $\al>\tfrac12$, we may fix $\gamma$ with
\begin{equation}\label{eq:gammawindow}
1-\al\ <\ \gamma\ <\ \tfrac12 .
\end{equation}
By Lemma \ref{lem:hbdd}\,\textup{(i)}, $\hat h$ is a continuous linear map from
$C^{0,\gamma}[0,T]$ to $C[0,T]$, and by Theorem \ref{theo:kac},
\[\sT\Longrightarrow W\] in $C^{0,\gamma}[0,T]$. The continuous mapping theorem
\citep[Theorem~2.7]{Billingsley} and \eqref{eq:XhT} yield
$X^{(\al)}-x\rme^{-\beta\,\cdot}\Longrightarrow \hat h[W]$ in $C[0,T]$.

\emph{Step 2\textup: identification of $\hat h[W]$.}
For $\zeta>0$ let $K^{(\zeta)}(u):=K_\al(u\vee\zeta)$, a Lipschitz kernel on
$[0,T]$, and let $\hat h^{(\zeta)}$ be the operator \eqref{eq:byparts} built from
it. Since $K^{(\zeta)}=K_\al$ on $[\zeta,\infty)$ and $K^{(\zeta)}$ is constant
on $[0,\zeta]$, for every $g\in C^\gamma[0,T]$
\[
\bigl\|\hat h^{(\zeta)}[g]-\hat h[g]\bigr\|_\infty
\le[g]_\gamma\Bigl(2\sup_{0<u\le\zeta}|K_\al(u)|\,u^\gamma
+\int_0^{\zeta}|K_\al'(u)|\,u^\gamma\rmd u\Bigr)\ \longrightarrow\ 0
\qquad(\zeta\downarrow0),
\]
by \eqref{eq:C1C2}.
On the other hand, $K^{(\zeta)}$ being Lipschitz, the Wiener integral may be
integrated by parts\textup: for each $t$, almost surely,
\[
\int_0^tK^{(\zeta)}(t-s)\rmd W(s)
=K^{(\zeta)}(0)\,W(t)+\int_0^t\bigl(K^{(\zeta)}\bigr)'(t-s)\,W(s)\rmd s
=\hat h^{(\zeta)}[W](t),
\]
where the second equality follows by expanding \eqref{eq:byparts} and using
$W(0)=0$ and
$\int_0^t(K^{(\zeta)})'(t-s)\rmd s=K^{(\zeta)}(t)-K^{(\zeta)}(0)$.

Now let $\zeta\downarrow0$. Almost surely $W\in C^\gamma[0,T]$, so
$\hat h^{(\zeta)}[W](t)\to \hat h[W](t)$ almost surely. Moreover
\[
\bigl\|K^{(\zeta)}-K_\al\bigr\|^2_{L^2(0,t)}
=\int_0^{\zeta\wedge t}\bigl(K_\al(\zeta)-K_\al(u)\bigr)^2\rmd u
\ \longrightarrow\ 0,
\]
because $K_\al\in L^2(0,T)$ for $\al>\tfrac12$.
 Hence, by the It\^o isometry,
$\int_0^tK^{(\zeta)}(t-s)\rmd W(s)\to\int_0^tK_\al(t-s)\rmd W(s)$ in $L^2$.

The two limits coincide, so that
$\hat h[W](t)=\int_0^tK_\al(t-s)\rmd W(s)$ almost surely, for each $t$.
Since $\hat h[W]$ is continuous by Lemma \ref{lem:hbdd}\,\textup{(i)}, it is the
continuous version of the stochastic convolution, and \eqref{eq:kac-limit}
follows.

\emph{Step 3\textup: the limit process.}
A Wiener integral of a deterministic kernel is centred Gaussian, and the
It\^o isometry gives the covariance. 

For $\beta=0$ one has
$K_\al(u)=u^{\al-1}/\Gamma(\al)$ and
$\int_0^tK_\al^2=t^{2\al-1}/\bigl((2\al-1)\Gamma(\al)^2\bigr)$.
\end{proof}

The condition $\al>\tfrac12$ has been used twice\textup: for the window
\eqref{eq:gammawindow} to be non-empty, and for $K_\al$ to be
square-integrable, so that the Wiener integral in \eqref{eq:kac-limit} exists.

\subsection{The regime $\al\le\tfrac12$}\label{sec:kac-small}

For $\al\le\tfrac12$ the kernel $K_\al$ is not square-integrable at the origin,
and the Gaussian process in \eqref{eq:kac-limit} is not defined. The next
proposition shows that, under \eqref{def:scaling} and without a further
normalization, $X^{(\al)}$ has no limit in law, not even along a subsequence.

Two facts are used in its proof. By \eqref{eq:sigcorr} with $m=1$, for every
integrable $g$ on $[0,t]$ and whatever the initial law,
\begin{equation}\label{eq:kac-L2}
\EE\Bigl[\Bigl(\int_0^tg(s)\,\eta(s)\rmd s\Bigr)^2\Bigr]
=\int_{[0,t]^2}g(u)\,g(v)\,\rme^{-2\la|u-v|}\rmd u\,\rmd v .
\end{equation}
Furthermore, since by \eqref{eq:Kode} $K_\al(u)=u^{\al-1}/\Gamma(\al)-\beta K_{1+\al}(u)$, the
representation \eqref{eq:repr} takes the form
\begin{equation}\label{eq:XFR}
X^{(\al)}(t)=x\rme^{-\beta t}+\sF^{(\al)}(t)-\beta R(t),
\end{equation}
where $ R(t):=c\int_0^tK_{1+\al}(t-s)\,\eta(s)\rmd s .$

The function $K_{1+\al}$ is continuous on $[0,T]$, hence by \eqref{eq:kac-L2},
\begin{equation}\label{eq:Rbound}
\EE\bigl[R(t)^2\bigr]\le\|K_{1+\al}\|_{\infty,T}^2\,c^2\int_{[0,t]^2}\rme^{-2\la|u-v|}\rmd u\,\rmd v
\le\|K_{1+\al}\|_{\infty,T}^2\,T\,\frac{c^2}{\la},\qquad 0\le t\le T,
\end{equation}
which is bounded under \eqref{def:scaling}. So everything is decided by
$\sF^{(\al)}$.

\begin{proposition}[No scaling limit for $\al\le\tfrac12$]\label{prop:nolimit}
Let $0<\al\le\tfrac12$ and $0<t\le T$. Under \eqref{def:scaling}\textup,
\begin{equation}\label{eq:F2}
\EE\bigl[\sF^{(\al)}(t)^2\bigr]\ \sim\
\begin{cases}
\dfrac{(2\la)^{1-2\al}}{2\cos(\pi\al)}, & 0<\al<\tfrac12,\\[2ex]
\dfrac{\log(2\la t)}{\pi}, & \al=\tfrac12,
\end{cases}
\end{equation}
and the laws of $X^{(\al)}(t)$ are not tight. In particular\textup,
$X^{(\al)}$ has no limit in law in $C[0,T]$\textup, not even along a
subsequence.
\end{proposition}

\begin{proof}
\emph{Step 1\textup: the second moment.} 
Since
$\sF^{(\al)}(t)=\int_0^tg(s)\eta(s)\rmd s$
with $g(s)=c\,(t-s)^{\al-1}/\Gamma(\al)$,
\eqref{eq:kac-L2} and the substitution $y=2\la(t-u)$,
$z=2\la(t-v)$ give
\[
\EE\bigl[\sF^{(\al)}(t)^2\bigr]=\frac{c^2\,(2\la)^{-2\al}}{\Gamma(\al)^2}\,A_{2\la t}(\al),
\qquad
A_U(\al):=\int_{[0,U]^2}(yz)^{\al-1}\rme^{-|y-z|}\rmd y\,\rmd z .
\]
Let $\al<\tfrac12$. Then $A_U(\al)$ increases to $A_\infty(\al)$\textup, and with
$w=|y-z|$ and $s=(y\wedge z)/w$, so that $yz=s(1+s)w^2$ and the Jacobian is
$w$, we obtain
\[
A_\infty(\al)=2\int_0^\infty\rme^{-w}w^{2\al-1}\rmd w\int_0^\infty s^{\al-1}(1+s)^{\al-1}\rmd s
=2\,\Gamma(2\al)\,B(\al,1-2\al)=\frac{\Gamma(\al)^2}{\cos(\pi\al)};
\]
the last two steps use \citet[8.380(3)]{GR} and the reflection formula,
applied to $\Gamma(2\al)\Gamma(1-2\al)$ and to $\Gamma(\al)\Gamma(1-\al)$.
Since $A_{2\la t}(\al)\uparrow A_\infty(\al)$ and $c^2/\la\to1$, this gives
the first case of \eqref{eq:F2}.

Let $\al=\tfrac12$. 
The same substitution gives
$A_U(\tfrac12)=2\int_0^U\rme^{-w}\left(\int_w^U\bigl(y(y-w)\bigr)^{-1/2}\rmd y\right)\rmd w$, and
\[
\int_w^U\frac{\rmd y}{\sqrt{y(y-w)}}
=2\log\Bigl(\sqrt{U/w}+\sqrt{U/w-1}\Bigr)\in\bigl[\log(U/w),\,\log(4U/w)\bigr],
\]
so that $A_U(\tfrac12)=2\log U+O(1)$ as $U\to\infty$. Since $c^2/\la\to1$ and
$\Gamma(\tfrac12)^2=\pi$, this gives the second case of \eqref{eq:F2}.

\emph{Step 2\textup: the fourth moment.} We claim that
\begin{equation}\label{eq:F4}
\EE\bigl[\sF^{(\al)}(t)^4\bigr]\le3\,\Bigl(\EE\bigl[\sF^{(\al)}(t)^2\bigr]\Bigr)^2 .
\end{equation}
Let $u\in[0,t]^4$ have distinct coordinates $u_{(1)}<u_{(2)}<u_{(3)}<u_{(4)}$ in
increasing order. By \eqref{eq:sigcorr} with $m=2$,
$\EE\bigl[\prod_{i=1}^4\eta(u_i)\bigr]
=\rme^{-2\la(u_{(2)}-u_{(1)})}\,\rme^{-2\la(u_{(4)}-u_{(3)})}$. This is one of the
three terms of
\[
\rme^{-2\la(|u_1-u_2|+|u_3-u_4|)}+\rme^{-2\la(|u_1-u_3|+|u_2-u_4|)}
+\rme^{-2\la(|u_1-u_4|+|u_2-u_3|)},
\]
and the other two are positive. Multiplying by
$\prod_{i}c\,(t-u_i)^{\al-1}/\Gamma(\al)>0$ and integrating over $[0,t]^4$, the
left-hand side becomes $\EE[\sF^{(\al)}(t)^4]$ and, by \eqref{eq:kac-L2}, each of
the three terms becomes $\bigl(\EE[\sF^{(\al)}(t)^2]\bigr)^2$.

\emph{Step 3\textup: the laws are not tight.} The Paley--Zygmund inequality
$\PP\{V\ge\tfrac12\EE V\}\ge\tfrac14(\EE V)^2/\EE[V^2]$, applied to
$V=\sF^{(\al)}(t)^2$, gives with \eqref{eq:F4}
\[
\PP\Bigl\{\bigl|\sF^{(\al)}(t)\bigr|\ge\bigl(\tfrac12\EE[\sF^{(\al)}(t)^2]\bigr)^{1/2}\Bigr\}\ \ge\ \frac1{12}.
\]
By \eqref{eq:Rbound} and Chebyshev's inequality there is $L<\infty$ such that
$\PP\{|x\rme^{-\beta t}-\beta R(t)|>L\}\le\tfrac1{24}$ for all large $\la$. Hence,
by \eqref{eq:XFR},
\[
\PP\Bigl\{\bigl|X^{(\al)}(t)\bigr|\ge\bigl(\tfrac12\EE[\sF^{(\al)}(t)^2]\bigr)^{1/2}-L\Bigr\}\ \ge\ \frac1{24},
\]
and the threshold tends to infinity by \eqref{eq:F2}. The same holds along every
subsequence, and convergence in law in $C[0,T]$ would make the laws of
$X^{(\al)}(t)$ tight.
\end{proof}

For $\al<\tfrac12$ the proposition also gives the size of $X^{(\al)}(t)$, which
is of order $\la^{1/2-\al}$\textup; at $\al=\tfrac12$ the divergence is only
logarithmic.

\section{Concluding remarks}\label{sec:conclusions}
\setcounter{equation}{0}

The fractional Kac--Ornstein--Uhlenbeck process is represented by the same formula
as its Gaussian counterpart: the initial condition decays exponentially, and
the noise is filtered through a kernel. What changes is the kernel. For
$0<\al<1$ and $\beta>0$ the Kummer kernel $K_\al$ changes sign exactly once and
has zero total mass, and the results of this paper are consequences of that
single fact. It confines the process to an explicit band, which is exactly the
support of its law and which depends on the two velocities only through their
difference. It makes the two conditional means cross, and it makes the
stationary version antipersistent: the spectral density vanishes at the origin,
and the autocovariance has zero integral. None of this occurs at $\al=1$, where
the kernel is the exponential, or at $\beta=0$, where it is the
Riemann--Liouville kernel. The interior of the parameter square is not an
interpolation between its two edges.

The first and second moments can be expressed explicitly through one family of
Kummer kernels, and the stationary variance is an elementary function of the
parameters. The stationary mean $\kappa$ of the velocity, which the limiting band does not see, survives in the law of large
numbers: $t^{-\al}\int_0^tX^{(\al)}(s)\rmd s$ converges almost surely to
$\kappa/(\beta\,\Gamma(1+\al))$.

The value $\al=\tfrac12$ is a threshold in two respects. Under Kac's scaling
the symmetric process converges to the Gaussian process
$x\rme^{-\beta t}+\int_0^tK_\al(t-s)\rmd W(s)$ when $\al>\tfrac12$, and has no
limit in law when $\al\le\tfrac12$. Moreover, for $\tfrac12<\al<1$ the low-frequency
behaviour of the stationary process is that of a fractional Gaussian noise of
Hurst index $H=\al-\tfrac12$, which is also the index of the limit when
$\beta=0$. In that range the fKOU process is thus driven by a bounded velocity and
has paths of bounded variation, yet under Kac's scaling it behaves like a Gaussian process of Hurst index $\al-\tfrac12<\tfrac12$.

Two extensions are natural. The first is the law of $X^{(\al)}(t)$ inside the
band, of which only the support is determined here. The second is the scaling
limit for non-symmetric velocities. We are preparing a series of publications on
this topic.

\backmatter

\section*{Statements and Declarations}

\bmhead{Funding}
The work of NR was partially supported by the Russian Science Foundation
and the Chelyabinsk Region, grant 26-21-20002, https://rscf.ru/project/26-21-20002.

\bmhead{Competing interests}
The authors have no competing interests to declare that are relevant to the
content of this article.

\bmhead{Author contributions}
Nikita Ratanov proposed the original idea of the study, and Oscar
L\'opez-Alfonso developed it. Both authors then formalized and extended the
results, read and revised the manuscript, made changes to it, and approved the
final version.

\bmhead{Use of AI tools}
The first draft of the manuscript, containing the original idea and its first
development, was written by the authors without the use of AI tools. In the
subsequent stages the authors used Claude (Anthropic) to assist with
formalizing the arguments, checking calculations, verifying references against
the original sources, and drafting and editing parts of the text. The authors
reviewed and edited all content and take full responsibility for the content
of the publication.

\bmhead{Data availability}
No data were generated or analysed. Every figure in the paper is produced by
scripts which are available from the authors on request.

\bibliography{fKOU1-JSP-25-11-2026}


\begin{thebibliography}{29}
\ifx \bisbn   \undefined \def \bisbn  #1{ISBN #1}\fi
\ifx \binits  \undefined \def \binits#1{#1}\fi
\ifx \bauthor  \undefined \def \bauthor#1{#1}\fi
\ifx \batitle  \undefined \def \batitle#1{#1}\fi
\ifx \bjtitle  \undefined \def \bjtitle#1{#1}\fi
\ifx \bvolume  \undefined \def \bvolume#1{\textbf{#1}}\fi
\ifx \byear  \undefined \def \byear#1{#1}\fi
\ifx \bissue  \undefined \def \bissue#1{#1}\fi
\ifx \bfpage  \undefined \def \bfpage#1{#1}\fi
\ifx \blpage  \undefined \def \blpage #1{#1}\fi
\ifx \burl  \undefined \def \burl#1{\textsf{#1}}\fi
\ifx \doiurl  \undefined \def \doiurl#1{\url{https://doi.org/#1}}\fi
\ifx \betal  \undefined \def \betal{\textit{et al.}}\fi
\ifx \binstitute  \undefined \def \binstitute#1{#1}\fi
\ifx \binstitutionaled  \undefined \def \binstitutionaled#1{#1}\fi
\ifx \bctitle  \undefined \def \bctitle#1{#1}\fi
\ifx \beditor  \undefined \def \beditor#1{#1}\fi
\ifx \bpublisher  \undefined \def \bpublisher#1{#1}\fi
\ifx \bbtitle  \undefined \def \bbtitle#1{#1}\fi
\ifx \bedition  \undefined \def \bedition#1{#1}\fi
\ifx \bseriesno  \undefined \def \bseriesno#1{#1}\fi
\ifx \blocation  \undefined \def \blocation#1{#1}\fi
\ifx \bsertitle  \undefined \def \bsertitle#1{#1}\fi
\ifx \bsnm \undefined \def \bsnm#1{#1}\fi
\ifx \bsuffix \undefined \def \bsuffix#1{#1}\fi
\ifx \bparticle \undefined \def \bparticle#1{#1}\fi
\ifx \barticle \undefined \def \barticle#1{#1}\fi
\bibcommenthead
\ifx \bconfdate \undefined \def \bconfdate #1{#1}\fi
\ifx \botherref \undefined \def \botherref #1{#1}\fi
\ifx \url \undefined \def \url#1{\textsf{#1}}\fi
\ifx \bchapter \undefined \def \bchapter#1{#1}\fi
\ifx \bbook \undefined \def \bbook#1{#1}\fi
\ifx \bcomment \undefined \def \bcomment#1{#1}\fi
\ifx \oauthor \undefined \def \oauthor#1{#1}\fi
\ifx \citeauthoryear \undefined \def \citeauthoryear#1{#1}\fi
\ifx \endbibitem  \undefined \def \endbibitem {}\fi
\ifx \bconflocation  \undefined \def \bconflocation#1{#1}\fi
\ifx \arxivurl  \undefined \def \arxivurl#1{\textsf{#1}}\fi
\csname PreBibitemsHook\endcsname

\bibitem[\protect\citeauthoryear{Ratanov}{2021}]{KOU1}
\begin{barticle}
\bauthor{\bsnm{Ratanov}, \binits{N.}}:
\batitle{{O}rnstein--{U}hlenbeck processes of bounded variation}.
\bjtitle{Methodol. Comput. Appl. Probab.}
\bvolume{23}(\bissue{3}),
\bfpage{925}--\blpage{946}
(\byear{2021})
\doiurl{10.1007/s11009-020-09794-x}
\end{barticle}
\endbibitem

\bibitem[\protect\citeauthoryear{Ratanov}{2022}]{KOU2}
\begin{barticle}
\bauthor{\bsnm{Ratanov}, \binits{N.}}:
\batitle{{K}ac--{O}rnstein--{U}hlenbeck processes: stationary distributions and
  exponential functionals}.
\bjtitle{Methodol. Comput. Appl. Probab.}
\bvolume{24}(\bissue{4}),
\bfpage{2703}--\blpage{2721}
(\byear{2022})
\doiurl{10.1007/s11009-022-09956-z}
\end{barticle}
\endbibitem

\bibitem[\protect\citeauthoryear{Ratanov}{2025}]{FCAA}
\begin{botherref}
\oauthor{\bsnm{Ratanov}, \binits{N.}}:
{Note on a fractionally integrated telegraph process}.
SSRN preprint
(2025).
\doiurl{10.2139/ssrn.5192391}
\end{botherref}
\endbibitem

\bibitem[\protect\citeauthoryear{Haubold et~al.}{2011}]{saxena}
\begin{botherref}
\oauthor{\bsnm{Haubold}, \binits{H.J.}},
\oauthor{\bsnm{Mathai}, \binits{A.M.}},
\oauthor{\bsnm{Saxena}, \binits{R.K.}}:
{M}ittag-{L}effler functions and their applications.
J. Appl. Math.
\textbf{2011}
(2011)
\doiurl{10.1155/2011/298628} .
Art. ID 298628, 51 pp.
\end{botherref}
\endbibitem

\bibitem[\protect\citeauthoryear{Miller and Ross}{1993}]{MillerRoss}
\begin{bbook}
\bauthor{\bsnm{Miller}, \binits{K.S.}},
\bauthor{\bsnm{Ross}, \binits{B.}}:
\bbtitle{An Introduction to the Fractional Calculus and Fractional Differential
  Equations}.
\bpublisher{Wiley},
\blocation{New York}
(\byear{1993})
\end{bbook}
\endbibitem

\bibitem[\protect\citeauthoryear{Anh et~al.}{2019}]{anh-fract}
\begin{barticle}
\bauthor{\bsnm{Anh}, \binits{P.T.}},
\bauthor{\bsnm{Doan}, \binits{T.S.}},
\bauthor{\bsnm{Huong}, \binits{P.T.}}:
\batitle{A variation of constant formula for {C}aputo fractional stochastic
  differential equations}.
\bjtitle{Statist. Probab. Lett.}
\bvolume{145},
\bfpage{351}--\blpage{358}
(\byear{2019})
\doiurl{10.1016/j.spl.2018.10.010}
\end{barticle}
\endbibitem

\bibitem[\protect\citeauthoryear{Ascione et~al.}{2024}]{PirozziMishura}
\begin{bbook}
\bauthor{\bsnm{Ascione}, \binits{G.}},
\bauthor{\bsnm{Mishura}, \binits{Y.}},
\bauthor{\bsnm{Pirozzi}, \binits{E.}}:
\bbtitle{Fractional Deterministic and Stochastic Calculus}.
\bsertitle{De Gruyter Series in Probability and Stochastics}.
\bpublisher{De Gruyter},
\blocation{Berlin, Boston}
(\byear{2024}).
\doiurl{10.1515/9783110780017}
\end{bbook}
\endbibitem

\bibitem[\protect\citeauthoryear{Meerschaert and Straka}{2013}]{straka}
\begin{barticle}
\bauthor{\bsnm{Meerschaert}, \binits{M.M.}},
\bauthor{\bsnm{Straka}, \binits{P.}}:
\batitle{Inverse stable subordinators}.
\bjtitle{Math. Model. Nat. Phenom.}
\bvolume{8}(\bissue{2}),
\bfpage{1}--\blpage{16}
(\byear{2013})
\doiurl{10.1051/mmnp/20138201}
\end{barticle}
\endbibitem

\bibitem[\protect\citeauthoryear{Kondratiev}{2022}]{time-times}
\begin{barticle}
\bauthor{\bsnm{Kondratiev}, \binits{Y.}}:
\batitle{From time to times}.
\bjtitle{Fract. Calc. Appl. Anal.}
\bvolume{25}(\bissue{1}),
\bfpage{1}--\blpage{14}
(\byear{2022})
\doiurl{10.1007/s13540-021-00001-3}
\end{barticle}
\endbibitem

\bibitem[\protect\citeauthoryear{Leonenko and Pirozzi}{2022}]{LP21}
\begin{barticle}
\bauthor{\bsnm{Leonenko}, \binits{N.}},
\bauthor{\bsnm{Pirozzi}, \binits{E.}}:
\batitle{First passage times for some classes of fractional time-changed
  diffusions}.
\bjtitle{Stoch. Anal. Appl.}
\bvolume{40}(\bissue{4}),
\bfpage{735}--\blpage{763}
(\byear{2022})
\doiurl{10.1080/07362994.2021.1953386}
\end{barticle}
\endbibitem

\bibitem[\protect\citeauthoryear{Leonenko and Pirozzi}{2025}]{LP}
\begin{barticle}
\bauthor{\bsnm{Leonenko}, \binits{N.}},
\bauthor{\bsnm{Pirozzi}, \binits{E.}}:
\batitle{The time-changed stochastic approach and fractionally integrated
  processes to model the actin-myosin interaction and dwell times}.
\bjtitle{Math. Biosci. Eng.}
\bvolume{22}(\bissue{4}),
\bfpage{1019}--\blpage{1054}
(\byear{2025})
\doiurl{10.3934/mbe.2025037}
\end{barticle}
\endbibitem

\bibitem[\protect\citeauthoryear{Angelani et~al.}{2024}]{ANG}
\begin{barticle}
\bauthor{\bsnm{Angelani}, \binits{L.}},
\bauthor{\bsnm{De~Gregorio}, \binits{A.}},
\bauthor{\bsnm{Garra}, \binits{R.}},
\bauthor{\bsnm{Iafrate}, \binits{F.}}:
\batitle{{Anomalous Random Flights and Time-Fractional Run-and-Tumble
  Equations}}.
\bjtitle{J. Stat. Phys.}
\bvolume{191},
\bfpage{129}
(\byear{2024})
\doiurl{10.1007/s10955-024-03344-2}
\end{barticle}
\endbibitem

\bibitem[\protect\citeauthoryear{Kac}{1974}]{Kac}
\begin{barticle}
\bauthor{\bsnm{Kac}, \binits{M.}}:
\batitle{A stochastic model related to the telegrapher's equation}.
\bjtitle{Rocky Mountain J. Math.}
\bvolume{4}(\bissue{3}),
\bfpage{497}--\blpage{509}
(\byear{1974})
\doiurl{10.1216/RMJ-1974-4-3-497}
\end{barticle}
\endbibitem

\bibitem[\protect\citeauthoryear{Ratanov and Kolesnik}{2022}]{RK}
\begin{bbook}
\bauthor{\bsnm{Ratanov}, \binits{N.}},
\bauthor{\bsnm{Kolesnik}, \binits{A.D.}}:
\bbtitle{Telegraph Processes and Option Pricing},
\bedition{2nd} edn.
\bpublisher{Springer},
\blocation{Berlin, Heidelberg}
(\byear{2022}).
\doiurl{10.1007/978-3-662-65827-7}
\end{bbook}
\endbibitem

\bibitem[\protect\citeauthoryear{Ratanov}{1999}]{MPRF}
\begin{barticle}
\bauthor{\bsnm{Ratanov}, \binits{N.}}:
\batitle{Telegraph evolutions in inhomogeneous media}.
\bjtitle{Markov Process. Related Fields}
\bvolume{5}(\bissue{1}),
\bfpage{53}--\blpage{68}
(\byear{1999})
\end{barticle}
\endbibitem

\bibitem[\protect\citeauthoryear{Kelbert and Orsingher}{1994}]{Kelbert}
\begin{barticle}
\bauthor{\bsnm{Kelbert}, \binits{M.Y.}},
\bauthor{\bsnm{Orsingher}, \binits{E.}}:
\batitle{On a telegraph-type equation with non-constant coefficients emerging
  in randomly accelerated motions}.
\bjtitle{Problems Inform. Transmission}
\bvolume{30}(\bissue{2}),
\bfpage{177}--\blpage{182}
(\byear{1994}).
\bcomment{Translated from Probl. Peredachi Inf. 30(2) (1994), 99--103}
\end{barticle}
\endbibitem

\bibitem[\protect\citeauthoryear{Ratanov and Turov}{2025}]{GSA25}
\begin{barticle}
\bauthor{\bsnm{Ratanov}, \binits{N.}},
\bauthor{\bsnm{Turov}, \binits{M.}}:
\batitle{On subordinated telegraph processes}.
\bjtitle{Global Stoch. Anal.}
\bvolume{12}(\bissue{5}),
\bfpage{38}--\blpage{44}
(\byear{2025})
\end{barticle}
\endbibitem

\bibitem[\protect\citeauthoryear{Ratanov}{2026}]{APT}
\begin{barticle}
\bauthor{\bsnm{Ratanov}, \binits{N.}}:
\batitle{Metastability for telegraph processes in a double-well potential}.
\bjtitle{J. Appl. Probab.}
\bvolume{63}(\bissue{1}),
\bfpage{204}--\blpage{224}
(\byear{2026})
\doiurl{10.1017/jpr.2025.10027}
\end{barticle}
\endbibitem

\bibitem[\protect\citeauthoryear{Br\'emaud}{2020}]{Bremaud}
\begin{bbook}
\bauthor{\bsnm{Br\'emaud}, \binits{P.}}:
\bbtitle{Markov Chains: {G}ibbs Fields, {M}onte {C}arlo Simulation and Queues},
\bedition{2nd} edn.
\bsertitle{Texts in Applied Mathematics},
vol. \bseriesno{31}.
\bpublisher{Springer},
\blocation{Cham}
(\byear{2020}).
\doiurl{10.1007/978-3-030-45982-6}
\end{bbook}
\endbibitem

\bibitem[\protect\citeauthoryear{Abramowitz and Stegun}{1972}]{AS}
\begin{bbook}
\beditor{\bsnm{Abramowitz}, \binits{M.}},
\beditor{\bsnm{Stegun}, \binits{I.A.}} (eds.):
\bbtitle{Handbook of Mathematical Functions with Formulas, Graphs, and
  Mathematical Tables}.
\bsertitle{National Bureau of Standards Applied Mathematics Series},
vol. \bseriesno{55}.
\bpublisher{U.S. Government Printing Office},
\blocation{Washington, DC}
(\byear{1972})
\end{bbook}
\endbibitem

\bibitem[\protect\citeauthoryear{Prudnikov et~al.}{1992}]{PBM5}
\begin{bbook}
\bauthor{\bsnm{Prudnikov}, \binits{A.P.}},
\bauthor{\bsnm{Brychkov}, \binits{Y.A.}},
\bauthor{\bsnm{Marichev}, \binits{O.I.}}:
\bbtitle{Integrals and Series, Vol.~5: Inverse {L}aplace Transforms}.
\bpublisher{Gordon and Breach},
\blocation{New York}
(\byear{1992})
\end{bbook}
\endbibitem

\bibitem[\protect\citeauthoryear{Beran et~al.}{2013}]{LM}
\begin{bbook}
\bauthor{\bsnm{Beran}, \binits{J.}},
\bauthor{\bsnm{Feng}, \binits{Y.}},
\bauthor{\bsnm{Ghosh}, \binits{S.}},
\bauthor{\bsnm{Kulik}, \binits{R.}}:
\bbtitle{Long-Memory Processes: Probabilistic Properties and Statistical
  Methods}.
\bpublisher{Springer},
\blocation{Berlin, Heidelberg}
(\byear{2013}).
\doiurl{10.1007/978-3-642-35512-7}
\end{bbook}
\endbibitem

\bibitem[\protect\citeauthoryear{Mandelbrot and Van~Ness}{1968}]{MvN}
\begin{barticle}
\bauthor{\bsnm{Mandelbrot}, \binits{B.B.}},
\bauthor{\bsnm{Van~Ness}, \binits{J.W.}}:
\batitle{Fractional {B}rownian motions, fractional noises and applications}.
\bjtitle{SIAM Rev.}
\bvolume{10}(\bissue{4}),
\bfpage{422}--\blpage{437}
(\byear{1968})
\doiurl{10.1137/1010093}
\end{barticle}
\endbibitem

\bibitem[\protect\citeauthoryear{Lamperti}{1962}]{Lamperti}
\begin{barticle}
\bauthor{\bsnm{Lamperti}, \binits{J.}}:
\batitle{On convergence of stochastic processes}.
\bjtitle{Trans. Amer. Math. Soc.}
\bvolume{104}(\bissue{3}),
\bfpage{430}--\blpage{435}
(\byear{1962})
\doiurl{10.1090/S0002-9947-1962-0143245-1}
\end{barticle}
\endbibitem

\bibitem[\protect\citeauthoryear{Hamadouche}{2000}]{Hamadouche}
\begin{barticle}
\bauthor{\bsnm{Hamadouche}, \binits{D.}}:
\batitle{Invariance principles in {H}\"older spaces}.
\bjtitle{Portugal. Math.}
\bvolume{57}(\bissue{2}),
\bfpage{127}--\blpage{151}
(\byear{2000})
\end{barticle}
\endbibitem

\bibitem[\protect\citeauthoryear{Billingsley}{1999}]{Billingsley}
\begin{bbook}
\bauthor{\bsnm{Billingsley}, \binits{P.}}:
\bbtitle{Convergence of Probability Measures},
\bedition{2nd} edn.
\bsertitle{Wiley Series in Probability and Statistics}.
\bpublisher{Wiley},
\blocation{New York}
(\byear{1999}).
\doiurl{10.1002/9780470316962}
\end{bbook}
\endbibitem

\bibitem[\protect\citeauthoryear{L\'evy}{1953}]{Levy}
\begin{bbook}
\bauthor{\bsnm{L\'evy}, \binits{P.}}:
\bbtitle{Random Functions: General Theory with Special Reference to {L}aplacian
  Random Functions}.
\bsertitle{Univ. California Publications in Statistics},
vol. \bseriesno{1, no.~12}.
\bpublisher{Univ. California Press},
\blocation{Berkeley}
(\byear{1953})
\end{bbook}
\endbibitem

\bibitem[\protect\citeauthoryear{Gradshteyn and Ryzhik}{1980}]{GR}
\begin{bbook}
\bauthor{\bsnm{Gradshteyn}, \binits{I.S.}},
\bauthor{\bsnm{Ryzhik}, \binits{I.M.}}:
\bbtitle{Table of Integrals, Series, and Products},
\bedition{Corrected and enlarged} edn.
\bpublisher{Academic Press},
\blocation{New York}
(\byear{1980}).
\bcomment{Edited by A. Jeffrey}
\end{bbook}
\endbibitem

\end{thebibliography}

\end{document}